\documentclass[12pt]{amsart}

\usepackage{amsmath,amsfonts,amssymb}

\theoremstyle{plain}
\newtheorem{lemma}{Lemma}[section]
\newtheorem{theorem}[lemma]{Theorem}

\newtheorem{proposition}[lemma]{Proposition} 

\theoremstyle{definition}
\newtheorem{remark}[lemma]{Remark}
\newtheorem{example}[lemma]{Example}
\newtheorem{definition}[lemma]{Definition}

\numberwithin{equation}{section}
\DeclareMathOperator{\Var}{Var}

\DeclareMathOperator{\diam}{diam}
\DeclareMathOperator{\cov}{cov}

\DeclareMathOperator{\Int}{Int}
\DeclareMathOperator{\Ext}{Ext}
\DeclareMathOperator{\Prob}{Pr}
\DeclareMathOperator{\solid}{solid}
\DeclareMathOperator{\main}{main}
\DeclareMathOperator{\far}{far}
\DeclareMathOperator{\near}{near}

\begin{document}

\newcommand{\omB}{\overline{\mathcal{B}}}
\newcommand{\oG}{\overline{G}}
\newcommand{\oL}{\overline{\Lambda}}
\newcommand{\graph}{(\oL,\omB(\Lambda))}
\newcommand{\ZZ}{\mathbb{Z}^{2}}
\newcommand{\ZZd}{\mathbb{Z}^{d}}
\newcommand{\RR}{\mathbb{R}^{2}}
\newcommand{\RRd}{\mathbb{R}^{d}}
\newcommand{\lra}{\leftrightarrow}
\newcommand{\lrad}{\overset{}{\leftrightarrow}}
\newcommand{\ep}{\epsilon}
\newcommand{\omU}{\omega_{\mathcal{U}}}
\newcommand{\omS}{\omega_{\mathcal{S}}}
\newcommand{\omT}{\omega_{\mathcal{T}}}
\newcommand{\omSt}{\tilde{\omega}_{\mathcal{S}}}
\newcommand{\ome}{\omega_{e}}
\newcommand{\omet}{\tilde{\omega}_{e}}
\newcommand{\zeU}{\zeta_{\mathcal{U}}}
\newcommand{\zeS}{\zeta_{\mathcal{S}}}
\newcommand{\zeSt}{\tilde{\zeta}_{\mathcal{S}}}
\newcommand{\zte}{\zeta_{e}}
\newcommand{\zeet}{\tilde{\zeta}_{e}}
\newcommand{\Ue}{\mathcal{S} \backslash \{e\}}
\newcommand{\mA}{\mathcal{A}}
\newcommand{\mB}{\mathcal{B}}
\newcommand{\mC}{\mathcal{C}}
\newcommand{\mD}{\mathcal{D}}
\newcommand{\mE}{\mathcal{E}}
\newcommand{\mF}{\mathcal{F}}
\newcommand{\mG}{\mathcal{G}}
\newcommand{\mH}{\mathcal{H}}
\newcommand{\mI}{\mathcal{I}}
\newcommand{\mJ}{\mathcal{J}}
\newcommand{\mK}{\mathcal{K}}
\newcommand{\mN}{\mathcal{N}}
\newcommand{\mQ}{\mathcal{Q}}
\newcommand{\mR}{\mathcal{R}}
\newcommand{\mS}{\mathcal{S}}
\newcommand{\mT}{\mathcal{T}}
\newcommand{\mU}{\mathcal{U}}
\newcommand{\mV}{\mathcal{V}}
\newcommand{\mW}{\mathcal{W}}
\newcommand{\mX}{\mathcal{X}}
\newcommand{\mY}{\mathcal{Y}}
\newcommand{\mZ}{\mathcal{Z}}
\newcommand{\mkB}{\mathfrak{B}}
\newcommand{\mkC}{\mathfrak{C}}
\newcommand{\mkD}{\mathfrak{D}}
\newcommand{\mkL}{\mathfrak{L}}
\newcommand{\mkM}{\mathfrak{M}}
\newcommand{\mkN}{\mathfrak{N}}
\newcommand{\mkS}{\mathfrak{S}}
\newcommand{\mkR}{\mathfrak{R}}
\newcommand{\omkR}{\overline{\mathfrak{R}}}
\newcommand{\mkV}{\mathfrak{V}}
\newcommand{\bs}{\backslash}
\newcommand{\half}{\frac{1}{2}}

\title[Separated-Occurrence Inequalities]
{Separated-Occurrence Inequalities\\
for Dependent Percolation and Ising Models}
\author{Kenneth S. Alexander}
\address{Department of Mathematics DRB 155\\
University of Southern California\\
Los Angeles, CA  90089-1113 USA}
\email{alexandr@math.usc.edu}
\thanks{Research supported by NSF grants DMS-9802368 and DMS-0103790.}

\date{\today}

\begin{abstract}
Separated-occurrence inequalities are variants for dependent lattice models of the
van den Berg-Kesten inequality for independent models.  They take the form $P(A
\circ_r B) \leq (1 + ce^{-\ep r})P(A)P(B)$, where $A \circ_r B$ is the event that
$A$ and $B$ occur at separation $r$ in a configuration $\omega$, that is, there
exist two random sets of bonds or sites separated by at least distance $r$, one set
responsible for the occurrence of the event $A$ in $\omega$, the other for the
occurrence of
$B$.  We establish such inequalities for subcritical FK models, and for Ising
models which are at supercritical temperature or have an external field, with $A$
and $B$ increasing or decreasing events.
\end{abstract}
\maketitle

\section{Introduction and Preliminaries} \label{S:intro}
We begin with an informal description; full definitions will be given below.
The van den Berg-Kesten inequality \cite{vdBK}, generalized by Reimer in
\cite{Re}, is among the most powerful tools available for the study of
independent percolation.  This inequality deals with disjoint occurrence of
two events, which for bond percolation means loosely that in some
configuration, there are two disjoint sets of bonds, one responsible for the
occurrence of an event $A$ and the other for the occurrence of another event 
$B$.  These disjoint sets need not be deterministic---they 
may depend on the configuration.  For example if, for some $x$ and $y$, $A$ is the
event that there is a path of open bonds from $x$ to $y$, then any such path is a
set of bonds responsible for the occurrence of $A$.
Letting $A \circ B$
denote the event that $A$ and $B$ occur disjointly, the van den Berg-Kesten
inequality (specialized to the context of independent percolation) states that if
$A$ and $B$ are either increasing or decreasing, then 
\begin{equation} \label{E:vdBK1}
  P(A \circ B) \leq P(A)P(B),
\end{equation}
complementing the Harris-FKG inequality \cite{Ha} which states that if $A,B$ are
both increasing or both decreasing,
\[
  P(A \cap B) \geq P(A)P(B).
\]
Reimer's generalization extends (\ref{E:vdBK1}) to arbitrary $A$ and $B$.
For dependent percolation, such as the Fortuin-Kasteleyn random cluster model
(briefly, the \emph{FK model}), and for spin systems, (\ref{E:vdBK1}) cannot be
true in general, even for increasing events.  
For example, in the FK model with (in standard notation) $q >
1$, if $A$ is the event that some bond $e$ is open and $B$ is the event that some
other bond $f$ is open, then $A$ and $B$ by definition can only occur disjointly,
but they may be strictly positively correlated.  Grimmett \cite{Gr95b} proved a
version of (\ref{E:vdBK1}) for the FK model, but with two different probability
measures on the right side.

If a lattice model has good mixing properties, though, we may hope that
(\ref{E:vdBK1}) is approximately true if we require that $A$ and $B$ occur not
just disjointly but well-separated from each other.  Specifically, we seek
inequalities of the form
\[
  P(A \text{ and } B \text{ occur at separation } r \text{ or more }) \leq
    (1 + Ce^{-\lambda r})P(A)P(B),
\]
where $C, \lambda$ are constants not depending on $A,B$; the precise definition of
the above event will be given below.  We call such an inequality a
\emph{separated-occurrence inequality}.  Existing results in this direction either
require that the locations where $A$ and $B$ occur be deterministic \cite{Al98mix}
or restrict $A$ or $B$ to be a quite special type of event \cite{Al00wu}.  Our aim
here is mainly to prove extensions of (\ref{E:vdBK1}) which do not have these
restrictions.  

Turning to more formal definitions, let $J$ and $\Delta$ be finite sets and $A, B
\subset J^{\Delta}$.  For $\omega \in J^{\Delta}$ and $\Theta \subset \Delta$, we
say that
$A$ \emph{occurs on} $\Theta$ in the configuration $\omega$ if 
\[
  \omega^{\prime} \in J^{\Delta}, \omega^{\prime}_{x} = \omega_{x} 
    \text{ for all } x \in \Theta \text{ implies } \omega^{\prime} \in A.
\]
For a (possibly random) set $\Theta = \Theta(\omega) \subset \Delta$, we say that
$A$ \emph{occurs only on} $\Theta$ if $\omega \in A$ implies that $A$ occurs on
$\Theta(\omega)$ in $\omega$.
$A$ and $B$ are said to \emph{occur disjointly} in $\omega$ if there exist
disjoint $\Theta, \Gamma \subset \Delta$ such that $A$ occurs on $\Theta$ in
$\omega$ and $B$ occurs on $\Gamma$ in $\omega$.  The event that $A$ and $B$
occur disjointly is denoted $A \circ B$.  A linear ordering of $J$ induces the
coordinate-wise partial ordering on $J^{\Delta}$; we then say an event $A$ is
\emph{increasing} if $\omega \in A, \omega \leq \omega^{\prime}$ imply
$\omega^{\prime} \in A$.  $A$ is \emph{decreasing} if its complement $A^{c}$ is
increasing.  The van den Berg-Kesten inequality
\cite{vdBK} states that (\ref{E:vdBK1}) holds for every product measure $P$ on
$J^{\Delta}$ and all $A, B$ which are either increasing or
decreasing.  

We consider now analogous concepts suited to dependent percolation and
lattice random fields.  By a \emph{site} we mean an element of $\ZZd$; sites $x$
and $y$ are \emph{adjacent} if $|y - x| = 1$.  Here $| \cdot |$ denotes
the Euclidean norm. 
By a \emph{bond} we mean an unordered pair $\langle xy
\rangle$ of adjacent sites.  When
convenient we view a bond as a closed line segment in $\mathbb{R}^{d}$.  For $R
\subset \RRd$ we let $\mB(R) = \{ b \in \mB(\ZZd): b \subset R \}$, except that
for $\Lambda \subset \ZZd$ we let $\mB(\Lambda) = \{ \langle xy \rangle: x,y
\in \Lambda \}$; the context will prevent any ambiguity.  We also write
$\omB(\Lambda) = \{ \langle xy \rangle: x \in \Lambda \}$, and
\[
  V(\mR) = \{x \in \ZZd: \langle xy \rangle \in \mR \text{ for
    some } y \}.
\]
For $\mU,\mV \subset \mB(\ZZd)$ we say that $\mU$ \emph{abuts} $\mV$ if $\mU \cap
\mV = \phi$ but $V(\mU) \cap V(\mV) \neq \phi$.  A \emph{bond configuration} on a
set $\mR$ of bonds is an element $\omega \in
\{0,1\}^{\mR}$; when convenient we view $\omega$ as a subset of $\mR$ or
as a subgraph of $(V(\mR),\mR)$. 
A bond $e$ is \emph{open} in the
configuration $\omega$ if $\omega_{e} = 1$, and \emph{closed} if $\omega_{e} = 0$.
Given $\rho \in \{0,1\}^{\mR^{c}}$ we define $(\omega \rho)
= (\omega\rho)_{\mR}$ to be the bond
configuration on the full lattice which coincides with $\omega$ on $\mR$
and with $\rho$ on $\mR^{c}$.  

For $x \in \ZZd$ let $Q(x) = x + [-\tfrac{1}{2},\tfrac{1}{2}]^{d}$, and for
$\Lambda \subset \ZZd$ let $Q(\Lambda) = \cup_{x \in \Lambda} Q(x)$.
A \emph{dual plaquette} (or \emph{dual bond}, in two dimensions) is a face of a
cube $Q(x)$ for some $x \in \ZZd$.  A \emph{dual site} is a point $x +
(\tfrac{1}{2},..,\tfrac{1}{2})$ with $x \in \ZZd$.; the set of all dual sites is
denoted $(\ZZd)^*$.  Each dual plaquette  perpendicularly bisects a unique
bond $e$; we then denote the dual plaquette by $e^{}$.  
The dual plaquette $e^{}$ is
defined to be open precisely when $e$ is closed; in this way we obtain a dual
configuration $\omega^{}$ of dual plaquettes for
 each bond configuration $\omega$.  A
\emph{dual surface} (or \emph{dual circuit}, in two dimensions) is the
boundary of a set $Q(\Lambda)$ for some
$\Lambda \subset \ZZd$ for which $\mB(\Lambda)$ is connected.  A dual
surface is
\emph{open} if all its dual plaquettes are open.

A \emph{path} is a sequence $\gamma = 
(x_{0},\langle x_{0}x_{1} \rangle, x_{1},\ldots
x_{n-1},\langle x_{n-1}x_{n} \rangle, x_{n})$ of alternating 
sites and bonds.  The path $\gamma$ is called \emph{open} if all bonds in
$\gamma$ are open.  Let $x \lra y$ denote the event 
that there is an open path from $x$
to $y$.  
The \emph{cluster} of a set $\Theta \subset \ZZd$ in a configuration
$\omega$ is 
\[
  C(\Theta,\omega) = \{ x \in \ZZd: x \lra \Theta \text{ in } \omega \}.
\]
We write $C_{x}(\omega)$ for $C(\{x\},\omega)$, and when a bond $e$ is open in
$\omega$ we write $C_e(\omega)$ for $C(V(e),\omega)$.

For $\mR$ a set of bonds and $x,y \in V(\mR)$ we let $d_{\mR}(x,y)$ denote the
minimum length among all paths in $\mR$ from $x$ to $y$.  This determines a
distance between sets of sites, and for sets $\mE, \mF$ of bonds we define
$d_{\mR}(\mE,\mF) = d_{\mR}(V(\mE),V(\mF))$.  $\diam_{\mR}(\cdot)$ denotes
diameter for the distance $d_{\mR}$.  For $A, B \subset
\{0,1\}^{\mR}$ and $r > 0$ we say that $A$ and $B$ \emph{occur at
separation} $r$ in the bond configuration $\omega$ if there exist $\mathcal{E},
\mathcal{F} \subset \mR$ with $d_{\mR}(\mathcal{E},\mathcal{F}) \geq r$ such
that $A$ occurs on $\mathcal{E}$ in $\omega$ and $B$ occurs on $\mathcal{F}$ in
$\omega$.   
The event that $A$ and $B$ occur at separation $r$ is denoted $A \circ_{r} B$.

By a \emph{bond percolation model} we mean a probability measure $P$ on
$\{0,1\}^{\mR}$ for some $\mR \subset \mB(\ZZd)$.  When $\mR = \mB(\ZZd)$, the
conditional distributions for the model are denoted
\[
  P_{\mR,\rho} = P(\cdot \mid \omega_{e}
  = \rho_{e} \text{ for all } e \in \mR^{c}),
\]
where $\mR \subset \mathcal{B}(\ZZd)$.  
We write $\rho^{i}$ for the bond configuration consisting of all $i$'s, $i = 0,1$. 
When used as boundary conditions, 
$\rho^1$ and $\rho^0$ are called \emph{wired} and \emph{free} respectively,
and we sometimes write $P_{\mR,w}, P_{\mR,f}$ for $P_{\mR,\rho^1}, P_{\mR,\rho^0}$
respectively. Write $\omega_{\mathcal{D}}$ for $\{\omega_{e}: e \in \mathcal{D}\}$,
and let $\mathcal{G}_{\mathcal{D}}$ denote the 
$\sigma$-algebra generated by $\omega_{\mathcal{D}}$.

For $P$ a bond percolation model on $\mB(\ZZd)$,
we say that $P$ has \emph{bounded energy} if there exists $p_{0} > 0$ such that 
\begin{equation} \label{E:bddener}
  p_{0} < P \bigl( \omega_{e} = 1 \mid (\omega_{b}, b \neq e) \bigr) < 1 - p_{0}
    \qquad \text{for all } e \text{ and all } (\omega_{b}, b \neq e).
\end{equation}
We say that $P$ has \emph{exponential decay of connectivity} if there
exist $C, \lambda > 0$ such that
\[
  P(x \lra y) \leq C e^{-\lambda |y - x|} \quad \text{for all } x, y \in \ZZd.
\]
In two dimensions, $P$ has \emph{exponential decay of dual connectivity} if there
exist $C, \lambda > 0$ such that
\[
  P(x \lra y \text{ via a path of open dual bonds}) 
  \leq C e^{-\lambda |y - x|} \quad \text{for all } x, y \in (\ZZd)^*.
\]
$P$ has the \emph{weak mixing property} if for some $C, \lambda > 0$, for all
finite sets $\mathcal{D},\mathcal{E}$ with $\mathcal{D} \subset \mathcal{E}$,
\begin{align} 
  \sup \{\Var(&P_{\mathcal{E},\rho}(\omega_{\mathcal{D}} \in \cdot),
    P_{\mathcal{E},\rho^{\prime}}(\omega_{\mathcal{D}} \in \cdot)): 
    \rho, \rho^{\prime} \in \{0,1\}^{\mathcal{E}^{c}}\} \notag \\
  &\leq C \sum_{x \in V(\mathcal{D}),y \in V(\mathcal{E}^{c})} 
    e^{-\lambda |x - y|},
\notag
\end{align}
where $\Var(\cdot,\cdot)$ denotes total variation distance between measures.  
Roughly, weak mixing means that the influence of the boundary condition on a
finite region decays  exponentially with distance from that region. 
Equivalently, for some 
$C, \lambda > 0$, for all sets $\mathcal{E}, \mathcal{F} \subset 
\mathcal{B}(\ZZd)$,
\begin{align} \label{E:weakmix}
  \sup \{|&P(E \mid F) - P(E)|:  E \in \mathcal{G}_{\mathcal{E}}, F \in
    \mathcal{G}_{\mathcal{F}}, P(F) > 0\} \\
  &\leq C \sum_{x \in V(\mathcal{E}),y \in V(\mathcal{F})} 
    e^{-\lambda |x - y|}. \notag
\end{align}
$P$ has the \emph{ratio weak mixing property} if for some $C, \lambda > 0$,
for all sets $\mathcal{E}, \mathcal{F} \subset 
\mathcal{B}(\ZZd)$,
\begin{align} \label{E:rweakmix}
  \sup &\left\{ \left| \frac{P(E \cap F)}{P(E)P(F)} - 1 \right| : E \in 
    \mathcal{G}_{\mathcal{E}}, F \in
    \mathcal{G}_{\mathcal{F}}, P(E)P(F) > 0 \right\} \\
  &\leq C \sum_{x \in V(\mathcal{E}),y \in V(\mathcal{F})}
    e^{-\lambda |x - y|}, \notag
\end{align}
whenever the right side of (\ref{E:rweakmix}) is less than 1.
Note that (\ref{E:rweakmix}) is much stronger than (\ref{E:weakmix}) for $E,F$ for
which the probabilities on the left side of (\ref{E:weakmix}) are much smaller
than the right side of (\ref{E:weakmix}).  Also, the right side of
(\ref{E:weakmix}) or (\ref{E:rweakmix}) is small when $d(\mathcal{E},\mathcal{F})$
is a sufficiently large multiple of $\log \min(|\mathcal{E}|,|\mathcal{F}|)$.
Here $d(\cdot, \cdot)$ dentoes Euclidean distance.
It was shown in \cite{Al98mix} that for the FK model in two
dimensions, exponential decay of connectivity (in infinite volume, with wired
boundary) implies ratio weak mixing.

We can consider spin systems as well as percolation models, but because the
properties of increasing and decreasing events are central to our arguments, we
must restrict attention to systems in which the common spin space at each site is
(at least partially) ordered.  We will in fact consider only the most natural
example of this type:  the \emph{Ising model} on $\ZZd$, with single-spin space
$\{-1,1\}$ and  Hamiltonian
\[
  H_{\Lambda,\eta}(\sigma_{\Lambda}) = -\sum_{\langle xy \rangle \in
    \omB(\Lambda)} \delta_{[(\sigma\eta)_{\Lambda}(x)
    = (\sigma\eta)_{\Lambda}(y)]} - h\sum_{x \in \Lambda} \sigma_x
\]
for the model on $\Lambda$ with external field $h$ and boundary condition $\eta$.
Here for
site configurations $\sigma, \eta \in \{-1,1\}^{\ZZd}$ we write $\sigma_{\Lambda}$
for $(\sigma_{x}, x \in \Lambda)$, and $(\sigma\eta)_{\Lambda}$
for the configuration which
coincides with $\sigma$ on $\Lambda$ and with $\eta$ on $\Lambda^c$. 
The corresponding finite-volume Gibbs distribution at inverse temperature $\beta$
is given by
\[
  \mu_{\Lambda,\eta}^{\beta,h}(\sigma_{\Lambda}) =
    \frac{1}{Z_{\Lambda,\eta}^{\beta,h}} e^{-\beta
    H_{\Lambda,\eta}(\sigma_{\Lambda})},
\]
where $Z_{\Lambda,\eta}^{\beta,h}$ is the partition function.  
When $h=0$ we denote the critical inverse temperature of the model
by $\beta_{c}(d)$.  The above definitions given for
bond percolation models extend straightforwardly to the Ising model,
as do the definitions to come in this section; we formulate things here mainly
for bond models to avoid unnecessary repetition.  
Let $\mathcal{H}_{\Lambda}$ denote the
$\sigma$-algebra generated by
$\sigma_{\Lambda}$.

Throughout the paper, $c_{1}, c_{2},...$ and $\ep_1, \ep_2,...$ denote constants
which depend only on the infinite-volume model, or family of finite-volume models,
under consideration.  We reserve $\ep_i$ for constants that are ``sufficiently
small.''

Weak mixing for a spin system or bond percolation model has a
variety of useful consequences, particularly in two dimensions; see \cite{MOS}. 
It directly implies for a bond percolation model that for some constants
$c_{i}$, for $\mathcal{E}, \mathcal{F} \subset 
\mathcal{B}(\ZZd), d(\mathcal{E},\mathcal{F}) \geq
r \geq c_{1} \log \min(|\mathcal{E}|,|\mathcal{F}|)$ and $A
\in \mathcal{G}_{\mathcal{E}}, B \in \mathcal{G}_{\mathcal{F}}$,
\begin{equation} \label{E:fixedloc}
  P(A \circ_{r} B) = P(A \cap B) \leq P(A)P(B) +
    c_{2}e^{-\ep_{1}r} \min(P(A),P(B)).
\end{equation}
Ratio weak mixing, in contrast, directly implies the stronger statement that
under the same conditions,
\begin{equation} \label{E:fixedloc2}
  P(A \cap B) = P(A \circ_{r} B) \leq (1 + c_{3}e^{-\ep_{2}r})P(A)P(B).
\end{equation}
Inequality (\ref{E:fixedloc2}) has been applied in \cite{Al00wu} and
\cite{Al00sp} in the context of coarse-graining of interfaces and their FK-model
analogs in two dimensions, to obtain approximate independence of separated segments
of the interface.  But (\ref{E:fixedloc2}) suffers from two deficiencies---first,
it is a statement about the infinite-volume measure and 
does not immediately apply when
$P$ is a finite-volume measure under a boundary condition.  More important, it
requires that the locations $\mathcal{E}$ and $\mathcal{F}$ be deterministic. 
This is problematic, for example, when one event, say $A$, is the event that $x
\lra y$ for some $x, y$,
because one cannot say in advance where the path will be.  If $\mathcal{E}$ can
be random, by contrast, one can take $\mathcal{E}$ to be the path itself.  This
particular event $A$ occurs on the cluster of a fixed deterministic $x$, so
the situation is remedied for $d = 2$ by Lemma 3.2 of
\cite{Al00wu}.  This lemma states that if $P$ (on $\{0,1\}^{\mB(\ZZ)}$)
has the FKG property and exponential
decay of connectivity, and satisfies other mild assumptions, if $\mathcal{F}
\subset \mathcal{B}(\ZZ), B \in \mathcal{G}_{\mathcal{F}}, r \geq c_{4} \log
|\mathcal{F}|, x \in \ZZ$, and $A$ is an event which occurs only on the
cluster $C_{x}$,  then 
\[
  P(A \circ_{r} B) \leq (1 + c_{5}e^{-\ep_{3}r})P(A)P(B).
\]
But such an event $A$ is a very special type; we seek here separated-occurrence
inequalities covering general increasing and decreasing events.
In generalizing the inequality in (\ref{E:fixedloc2}) to a separated-occurrence
inequality, in which
$\mathcal{E}$ and $\mathcal{F}$ are random, the hypothesis $r \geq c_{6} \log
\min(|\mathcal{E}|,|\mathcal{F}|)$ is no longer appropriate, so one must ask, how
large should one require  the separation $r$ to be?  When
working in a finite region $\mR$ of the lattice, $r \geq c \log
|\mR|$ for some $c$ may be reasonable, in view of the preceding, but for infinite
$\mR$ there is no obvious choice.  A solution from \cite{Al00wu}, which we choose
here, is to restrict one of the events $A,B$ to occur on 
a random subset of some deterministic finite region $\mD \subset \mR$ and
allow the other event to occur anywhere on $\mR$, including on another part of
$\mathcal{D}$.  We then require roughly that $r \geq c \log \diam(\mD)$, where
$\diam(\cdot)$ denotes Euclidean diameter.  (Using
$\diam(\mD)$ instead of $|\mD|$ avoids problems with geometrically irregular
$\mD$.)

We will need a version of the
Markov property for open dual surfaces, adapted to finite volumes.  Specifically, a
\emph{blocking partition} of a set $\mR$ of bonds is an ordered partition
$(\mX,\mY,\mZ)$ of
$\mR$ such that every path from $\mX$ to $\mZ$ includes at least one bond of
$\mY$; we call $\mY$ a \emph{blocking set}.  A probability measure $P$ on
$\{0,1\}^{\mR}$ has the
\emph{Markov property for blocking sets} if for every blocking partition
$(\mX,\mY,\mZ)$ of
$\mR$, 
\begin{align} \label{E:blocking}
  &\text{the configuration $\omega_{\mX}$ and $\omega_{\mZ}$ are conditionally
    independent} \\
  &\qquad \text{given that all bonds in $\mY$ are closed.} \notag
\end{align}
(In infinite volume the
usual Markov property takes $\mY$ to be the set of bonds crossing some dual
surface; our definition is a natural analog in finite volumes.)  This Markov
property says roughly that there is no ``communication via the boundary'' from one
side of a blocking set to the other.  Note we may have $\mY = \phi$ if no component
of $\mR$ intersects both $\mX$ and $\mZ$. If for some fixed
$\mW \subset \mR$, (\ref{E:blocking}) is valid under the additional assumption
that either $\mX$ or $\mZ$ contains $\mW$, we say
that the probability measure has the
\emph{Markov property for sets blocking} $\mW$.  For $\mW \subset \mE \subset
\mR$, we say that $\mW$ is \emph{blockable in} $\mE$ \emph{under} $P$ if
$P(\omega_{\mE} \in \cdot \mid \omega_{\mR \bs \mE} = \rho_{\mR \bs \mE}^i)$ has
the Markov property for sets blocking $\mW$, for $i = 0,1$.  This says roughly
that if we view $\mR \bs \mE$ as part of the boundary for $\mE$, when the boundary
condition on this partial boundary is free or wired, there is no communication
via the partial boundary from $\mW$ to the other side of a barrier blocking
$\mW$ in $\mE$; see Lemma \ref{L:markov}.

\input{epsf}
\begin{figure} 
\epsfxsize=3.5in
\begin{center}
\leavevmode
\epsffile{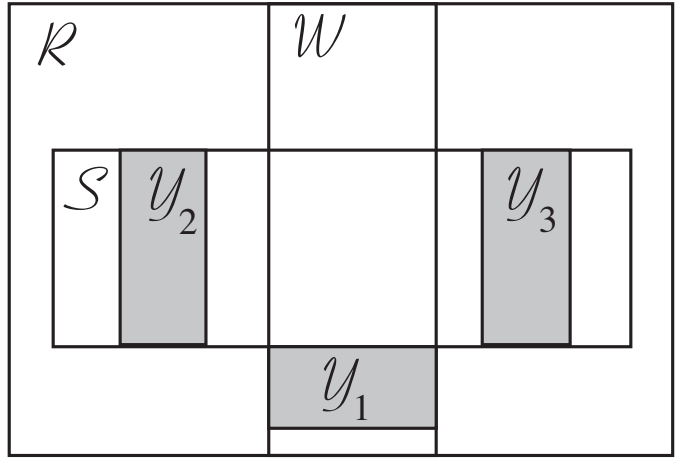}
\end{center}
\caption{A situation in which the Markov property holds for sets blocking
  $\mW$ but not for general blocking sets.  $\mS$ is the full horizontal
  rectangle, and $\mW$ is the full vertical rectangle.}
\label{F:markov}
\end{figure}

Figure \ref{F:markov} depicts a situation
for the FK model in which the Markov property holds for sets blocking $\mW$, but
not for general blocking sets.  
Let $\mE$ be the union of the horizontal rectangle $\mS$ and the vertical
rectangle $\mW$.  The boundary condition on $\mE$
consists of open bonds on the horseshoe-shaped regions comprising $\mR \bs (\mS
\cup \mW)$, and all bonds outside $\mR$ closed.   The Markov property does not
hold for the horizontal blocking set
$\mY_1$, because the presence or absense of a connection from one component of 
$\mR \bs (\mS \cup \mW)$ to the other component
underneath $\mY_1$ affects the probability of a  connection between these
components above
$\mY_1$.  However, the Markov property does hold for sets blocking
$\mE$, such as $\mY_2, \mY_3$ or $\mY_2 \cup \mY_3$; in fact $\mW$ is blockable in
$\mE$ under $\mW$.

For $\mR \subset \mB(\ZZd), r > 0$ and $x \in V(\mR)$ define the $d_{\mR}$-ball
\[
  B_{\mR}(x,r) = \{ b = \langle yz \rangle \in \mR: d_{\mR}(x,y) \leq r,
    d_{\mR}(x,z) \leq r \}.
\]
For $c > 1, r \geq 1$, a $c$-\emph{approximate r-neighborhood} of a
point $x \in \RRd$ in $\mR$ is a set
$\mN$ of bonds satisfying $B_{\mR}(x,r/c) \subset \mN \subset B_{\mR}(x,r)$.  A
$c$-approximate $r$-neighborhood of a bond $b$ is a $c$-approximate
$r$-neighborhood of an endpoint of $b$, say the one closest to the origin.  We make
analogous definitions for a set
$\Lambda$ of sites in place of the set $\mR$ of bonds.

A separated-occurrence inequality for finite volumes, in the form
(\ref{E:fixedloc2}), is most meaningful if the constants $c_{3},\ep_{2}$ are
uniform over some class of finite volumes and/or boundary conditions; for a single
finite-volume measure $P$ with bounded energy, there is always some choice of
$c_{3},\ep_{2}$ that makes the inequality trivially true.  As seen in the
related context of
\cite{Al01fin} for the FK model, properties like exponential decay and mixing may
hold in infinite volume but fail to hold (appropriately reformulated) uniformly in
finite volumes under certain kinds of boundary conditions.  Therefore it is
necessary to restrict to special classes of boundary conditions.  To give a
unified presentation without excessive numbers of cases, we use the following
formulation for bond models.  $\mkR$ is a collection of finite subsets of
$\mB(\ZZd)$, and for each
$\mR \in \mkR$ we have a collection $\mkM_{\mR}$ of probability measures on
$\{0,1\}^{\mR}$.  We seek
separated-occurrence inequalities which are uniform over
$\mkM = \cup_{\mR \in \mkR}\ \mkM_{\mR}$, that is, the constants $c_{3},\ep_{2}$ 
(as in (\ref{E:fixedloc2})) do not
depend on $P \in \mkM$.  Typically, each
$\mkM_{\mR}$ might consist of some class of bond or site boundary conditions, on
$\mR$, for some fixed bond percolation model on $\mB(\ZZd)$.
In addition we use auxiliary collections $\mkV_{\mR} \subset
\mkS_{\mR}$ of subsets of $\mR$ and an auxiliary collection $\mkN$ not depending
on $\mR$; we write $\mkV = \{ \mkV_{\mR}: \mR \in \mkR \}$
and $\mkS = \{ \mkS_{\mR}: \mR \in \mkR \}$.  For technical convenience we require
that $|\mV| \geq 2$ for all $\mV \in \mkV_{\mR}$.
We will refer to this setup---that is, to a designation of collections $\mkR,
\mkV, \mkS, \mkM, \mkN$ as described---as our \emph{standard bond percolation
setup}.  As we have mentioned, one of the two events in the separated-occurrence
inequality may be restricted to occur on a particular $\mD \subset \mR$; the sets
in $\mkV_{\mR}$ are neighborhoods of such $\mD$ in $\mR$, the sets in
$\mkS_{\mR}$ are subsets of such neighborhoods, and the sets in $\mkN$ are
approximate neighborhoods in sets $\mS \in \mkS_{\mR}$.  (One can take
$\mD = \mR$ and hence $\mkV_{\mR} = \{ \mR \}$ if desired, but allowing smaller
$\mD$ means we may reduce the required separation
$r$ between the events in question.)   

Like the proof in \cite{vdBF} of the van den Berg-Kesten inequality, our proofs of
separated-occurrence inequalities involve a process of splitting bonds one at a
time.  In our context the bonds to be split are those in some neighborhood $\mV
\in \mkV_{\mR}$ of $\mD$ in
$\mR$.  We may view this procedure as the filling of $\mV$ sequentially with
split bonds.  As this filling process proceeds we require, among other things, that
the  set $\mS \subset \mV$ of split bonds always satisfy
$\mS \in \mkS_{\mR}$.
We will need assumptions not just on the original
collections
$\mkR, \mkM$ but also on $\mkV,\mkS$ and certain augmented collections of measures
derived from $\mkM$. 
We have no general method for specifying the collections $\mkV, \mkS$ 
given a choice of original collections $\mkR$ and $\mkM$, but in the specific
examples we will consider, it is not particularly difficult to do so, as we will
see, though the methods are rather ad hoc.

\begin{definition} \label{D:enlarged}
Let $\mkR, \mkV, \mkS, \mkM$ be as above:  $\mkR$ is a collection of subsets of
$\mB(\ZZd)$, and for each $\mR \in \mkR,\ \mkV_{\mR} \subset \mkS_{\mR}$ are
collections of subsets of $\mR$ and $\mkM_{\mR}$ is a collection of 
probability measures on $\{0,1\}^{\mR}$.  
A collection $\mkN$ of subsets of $\mB(\ZZd)$ is a \emph{neighborhood collection
for} $\mkS$ if there exists $c>0$ such that for every $\mS \in \mkR \cup (\cup_{\mR
\in \mkR}\ 
\mkS_{\mR}), x \in V(\mS)$ and $r>0$, $\mkN$ contains a $c$-approximate
$r$-neighborhood of $x$ in $\mS$.  The corresponding \emph{augmented collections
of measures} are
\[
  \mkM_w^+(\mkN,\mkR) = \left\{ P(\omega_{\mN} \in \cdot \mid \omega_{\mR
    \bs \mN} = \rho_{\mR \bs \mN}^1): \mN \in \mkN, \mR \in \mkR, \mN \subset \mR,
    P \in \mkM \right\}
\]
and
\[
  \mkM^+(\mkS,\mkR) = \left\{ P(\omega_{\mS} \in \cdot \mid \omega_{\mR
    \bs \mS} = \rho_{\mR \bs \mS}^i): \mR \in \mkR, \mS \in \mkS_{\mR},
    P \in \mkM, i = 0,1 \right\}.
\]
Let $c>1$ and $r>0$; given $\mR \in \mkR$, $P \in
\mkM_{\mR}$, $\mS
\subset \mR$ and $e \in \mR \bs \mS$ we say that 
$e$ is \emph{neighborhood-appendable to} $\mS$ for $(c,r,\mR,P,\mkN)$
if there exists a $c$-approximate
$r$-neighborhood
$\mW \in \mkN$ of $e$ in $\mR$ such that either $\mW$ or $\mW \bs \mS$
is blockable in $\mS \cup \mW$ under $P$ (see Figure \ref{F:markov}.)
We say that $\mV \subset \mR$ is \emph{fillable compatibly with}
$(c,r,\mR,P,\mkS_{\mR},\mkN)$ if there exists an ordering $b_1,..,b_n$ of the
bonds of $\mV$ such that for all $0
\leq k < n$, we have (i) $\mS_k = \{b_1,..,b_k\} \in \mkS_{\mR}$, 
(ii) $b_{k+1}$ is neighborhood-appendable to $\mS_k$ for
$(c,r,\mR,P,\mkN)$, and (iii) 
every $e \in \mR \bs \mS$ is neighborhood-appendable to $\mS_k$ for
$(c,r/4,\mR,P,\mkN)$. For $R > 1$ we say that $(\mkR,\mkV,\mkS,\mkM,\mkN)$ is
\emph{filling-compatible at scale} $R$ if for some $c > 1$, for every $\mR \in
\mkR, \mV \in \mkV_{\mR}, P
\in \mkM_{\mR}$ and $1 \leq r < R$, $\mV$ is fillable compatibly with
$(c,r,\mR,P,\mkS_{\mR},\mkN)$.  
In this context we refer to such $c$ as \emph{suitable}.  We omit the wording ``at
scale $R$'' when $R = \infty$.
\end{definition}

Let us consider such properties for some natural classes $\mkR, \mkV, \mkS, \mkM$.
It is worth pointing out that we would certainly like the class $\mkR$ to be as 
large as possible, but it is less important that the classes $\mkV_{\mR}$ be
large, since the sets $\mV \in \mkV_{\mR}$ are used only to in effect enclose a
neighborhood of the set
$\mD$, where one of the events $A,B$ occurs, in a ``nice'' set that can serve as
the set of split bonds.  

For $\mR \subset \mB(\ZZd)$ and $k
\leq d$ let $\mR_k$ be the union of all
$k$-dimensional unit cubes having all edges in $\mR$, and let $\mR^{\solid} = 
\cup_{k \leq d} \mR_k$.  
We say that $\mR$ is \emph{simply lattice-connected}, abbreviated
\emph{SLC}, if
$\mR^{\solid}$ is simply connected.  For $d=2$, it is easy to see that $\mR$ is
SLC if and only if both
$\mR$ and $(\mR^c)^*$ are connected sets of bonds and dual bonds respectively,
which was the definition used in \cite{Al01fin}.  For $\Lambda \subset \ZZd$ we say
$\Lambda$ is simply lattice-connected if $\mB(\Lambda)$ is simply
lattice-connected.    For $d=2$ simple lattice-connectedness is a natural
assumption  in the context of uniform exponential decay and strong mixing, in view
of \cite{Al01fin}.   The SLC property also fits well with our need to form
enlarged collections, as Example \ref{Ex:slcprops} below shows.
For $\Gamma$ a self-avoiding lattice circuit let $\Int(\Gamma)$ denote the set of
all bonds strictly inside $\Gamma$ (excluding endpoints).
For $d=2$, we say that $\mR$ is \emph{circuit-bounded} if there exists a
self-avoiding lattice circuit, denoted $\Gamma_{\mR}$, such that $\mR =
\Gamma_{\mR} \cup
\Int(\Gamma_{\mR})$.  For $\mR$ circuit-bounded and $\mB \subset \Int(\mR)$, the
\emph{boundary closure} of $\mB$ in $\mR$ is $\mB \cup \{ \langle xy \rangle \in
\Gamma_{\mR}: x,y \in V(\mB) \}$.
We say that $\mR$ is \emph{lattice-convex} if it has the form $\{ b \in
\mB(\ZZd): b \subset C \}$ for some convex $C$.  If $C$ is a rectangle (which we
may assume has vertices in $\ZZd$) then we call $\mR$ a \emph{lattice rectangle},
and we let $\Gamma_{\mR}$ denote the outer surface $\{ b \in \mR: b \subset
\partial C \}$ and let $\Int(\mR) = \mR \bs \Gamma_{\mR}$.
If $\mR$ is a lattice rectangle (or a circuit-bounded set, or a $d_{\mA}$-ball
for some $\mA$),
$\mB \subset \mR$ and every bond of
$\mR \bs \mB$ abuts $\mR^c$, then we call $\mB$ an
\emph{approximate lattice rectangle} (or \emph{approximate circuit-bounded set},
or \emph{approximate} $d_{\mA}$-\emph{ball}.) and we refer to $\mR$ as the
\emph{completed} lattice rectangle (or circuit-bounded set, or $d_{\mA}$-ball.)  
A lattice rectangle or a circuit-bounded set of bonds $\mR$ is \emph{minimally fat}
if $\Int(\Gamma_{\mR})$ is connected.  An approximate lattice rectangle $\mR$, with
completed lattice rectangle $\tilde{\mR}$, is \emph{regular} if for every bond $e
\in \mR \cap \Gamma_{\tilde{\mR}}$, all bonds in $\Int(\tilde{\mR})$ which abut
$e$ (note there are at most two) are in $\mR$.

\begin{example} \label{Ex:slcprops}
Let $\mkR$ be the class of all minimally fat lattice rectangles in
$\mB(\ZZd)$, 
let $\mkS_{\mR}$ be the class of all approximate lattice rectangles $\mS \subset
\mR$ with $\mR \bs \mS$ connected,  and let $\mkV_{\mR}$ be the class of all
lattice rectangles $\mV
\subset \mR$ such that $\mR \bs \mV$ is
connected and $|\mV| \geq 2$.  (Note that for $d=2$ this excludes those $\mV$ which
intersect
$\Gamma_{\mR}$ only in two opposite faces.)  Let $\mkN$ be the class of all
approximate lattice rectangles in
$\mB(\ZZd)$.  Let
$P$ be an FK model on
$\mB(\ZZ)$ (see Section 2 for a description), and let $\mkM$ be the class of
all the associated finite-volume FK measures
$P_{\mR,*}$ with $\mR\in \mkR$ and
with $*=w$ or $f$.  Assume that either $i=1$ or there are no external fields.
In view of Lemma \ref{L:markov} below,
given $\mR \in \mkR,\mV \in \mkV_{\mR}$ and $r \geq 1$, we want to
fill $\mV$ in such a way that for all $k$, for $\mS_k$ as in Definition
\ref{D:enlarged} and $e \in \mR \bs \mS_k$, for some $d$-approximate
$r$-neighborhood $\mW$ of $e$ in
$\mR$, we have the following: (a) $\mS_k$ is an approximate lattice
rectangle, (b)
$\mR \bs \mS_k$ is connected, and (c) each component of $\mS_k
\bs \mW$ is an approximate latttice rectangle which abuts at most one component of
$\mR \bs (\mS_k \cup
\mW)$.  We can always
choose $\mW$ to be a lattice rectangle which either does not intersect $\mS_k$ or
intersects the interior of the completed $\mS_k$; this means that within (c) there
are two cases:  (c$^{\prime}$) $\mR \bs (\mS_k \cup
\mW)$ is connected and (c$^{\prime\prime}$) $\mR \bs (\mS_k \cup
\mW)$ has two components ($\mW$ ``slices through'' $\mR$.)
It follows from Lemma
\ref{L:markov} that under (c$^{\prime}$), $\mW \bs
\mS_k$ is blockable in $\mS_k \cup \mW$ under $P_{\mR,*}$, while under 
(c$^{\prime\prime}$),
$\mW$  is blockable in $\mS_k \cup \mW$ under $P_{\mR,*}$.  (Here $*=w$ or $f$.) 
It will follow that $(\mkR,\mkV,\mkS,\mkM,\mkN)$ is filling-compatible.  

In fact, to obtain (a), (b), (c), we can fill $\mV$ by first filling $\mV
\cap\Int(\mR)$ and then filling $\mV \cap \Gamma_{\mR}$.  It is not hard to see
that during the filling of
$\mV \cap\Int(\mR)$ we can keep
$\mS_k$ an SLC regular approximate lattice rectangle so that all of $\mR \bs
\mS_k$ is connected to
$\Gamma_{\mR}$, which ensures that $\mR \bs \mS_k$ is connected.  Since $\mW$ is
a lattice rectangle,
$\mR \bs \mW$ has at most 2 components $\mC_i$, and $\mC_i \cap \Gamma_{\mR}$ is
nonempty and connected for each $i$.  If there are 2 components $\mC_i$, then,
since $\mS_k$ is regular, each
$\mS_k \cap \mC_i$ is an approximate lattice rectangle and (c$^{\prime\prime}$)
holds.  If there is only one component $\mC_1 = \mR \bs \mW$, then, since $\mS_k
\subset \Int(\mR)$, we have
$\mR \bs (\mS_k \cup \mW)$ connected, i.e. (c$^{\prime}$) holds.  (Note that
for $d \geq 3$, when there is only one component $\mC_1$, we could have $\mS_k \bs
\mW$ connected but not SLC, meaning (c) fails.)  Thus all of our conditions (a),
(b) and (c) remain satisfied during the filling of $\mV \cap\Int(\mR)$.  If
$\mV \cap
\Gamma_{\mR}$ is empty we are done; if $\mV \cap \Gamma_{\mR}$ is nonempty then by
assumption we have $\mR \bs\mV$ connected.  But then it is not hard to see that we
can fill the rest of $\mV$ (i.e. fill $\mV \cap \Gamma_{\mR}$) keeping
$\Gamma_{\mR}
\cap \mS_k$ and $\Gamma_{\mR} \bs \mS_k$ both connected, with $\mS_k$ having an
approximately rectangular intersection with each face of $\mR$.  Following this
procedure we see that our preceding verification of (a), (b) and (c) remains
valid as we fill $\mV \cap \Gamma_{\mR}$.
\end{example}

It is ``straightforward but tedious'' to extend Example \ref{Ex:slcprops} to allow
general lattice-convex sets in $\mkR$ and $\mkV_{\mR}$, instead of just
rectangles.  In two dimensions we can be much more general, as the next example
shows.

\begin{example} \label{Ex:circbdd}
Let $\mkR$ be class of all circuit-bounded subsets of $\mB(\ZZ)$.
For $\mR \in \mkR$ let 
\begin{align}
  \mkV_{\mR} = \{ \mV \subset \mR:  \mV &\text{ is the boundary closure of}
    \notag \\
  &\text{some }
    d_{\Int(\mR)}-\text{ball in } \mR, \mR \bs \mV \text{ is 
    connected, and } |\mV| \geq 2 \}, \notag
\end{align}
and let
$\mkS_{\mR}$ be the class of all SLC $\mS \subset \mR$ with $\mR \bs \mS$
connected.  Let $\mkM$ be the class of
all the associated finite-volume FK measures
$P_{\mR,*}$ with $\mR\in \mkR$ and
with $*=w$ or $f$, and let $\mkN$ be the class of all finite SLC subsets of
$\mB(\ZZ)$.  Similarly to Example
\ref{Ex:slcprops}, given $\mR \in \mkR$ and $\mV = \{ b_1,..,b_n \}
\in \mkV_{\mR}$ we want to fill
$\mV$ so that (a) $\mS_k = \{ b_1,..,b_k \}$ is SLC, (b) $\mR \bs \mS_k$ is
connected, and (c) for $e \in \mR \bs \mS_k$, for some 
$d_{\mR}$-ball $\mW$ centered at an endpoint of $e$, each
component of $\mS_k \bs \mW$ abuts at most one component of $\mR \bs (\mS_k \cup
\mW)$. 

A useful observation about $d_{\mR}$-balls in the plane (also valid for
$d_{\Int(\mR)}$-balls) is as follows.  The ``outer surface inside $\Gamma_{\mR}$''
for a ball $B_{\mR}(z,k)$, by which we mean $\{ w \in
V(\mR) \bs V(\Gamma_{\mR}): d_{\mR}(z,w) = k \}$, must fall
along diagonal lines, i.e. it cannot contain two adjacent sites, as is easily seen.

We first fill $\mV \cap
\Int(\mR)$.  We order the bonds of $\mV \cap \Int(\mR) = \{ b_1,..,b_m \}$ in order
of increasing $d_{\Int(\mR)}$-distance from $x$, breaking ties in such a way that
$\mS_k =
\{ b_1,..,b_k \}$ remains SLC.  This means that each $\mS_k$ is an approximate
$d_{\Int(\mR)}$-ball. Now (a) is clear, and (b) follows from the fact that
$\Gamma_{\mR} \subset \mR \bs \mS_k$.  For (c), let $\mW =
B_{\mR}(y,l)$ for an endpoint $y$ of $e$ and some $l > 1$, and suppose
$\mS_k$ is an approximate
$d_{\Int(\mR)}$-ball, for which the completed $d_{\Int(\mR)}$-ball is
$B_{\Int(\mR)}(x,m)$ for some
$m>1$.  If $\mW \cap \mS_k = \phi$ then after shrinking $\mW$ slightly (say,
decrease $l$ by 2), $\mR \bs (\mS_k \cup \mW)$ becomes connected, and then
(c) is trivial.  Hence we assume $\mW \cap
\mS_k \neq \phi$.  If we enlarge $\mW$ slightly (say, increase $l$ by 2),
this ensures that $\mS_k \cup
\mW$ is connected, which means that each component of $\mR \bs (\mS_k \cup \mW)$
contains at most one segment of
$\Gamma_{\mR}$.  
We claim that every component of $\mR \bs (\mS_k
\cup \mW)$ contains exactly one segment
of $\gamma$.  Suppose instead that some component 
$\mA$ of $\mR \bs (\mS_k \cup \mW)$ does not intersect $\Gamma_{\mR}$; since
$\mS_k$ and $\mW$ are SLC, loosely speaking we conclude 
that $\mS_k$ and $\mW$ must
intersect on two sides of $\mA$ in such a way that their union surrounds $\mA$. 
More precisely, there must exist sites $u,v \in V(\mS_k) \cap V(\mW)$ and lattice
paths
$\gamma_{x \lra u}, \gamma_{x \lra v}, \gamma_{y \lra u}, \gamma_{y \lra v}$ of
minimal length (at most $m,m,l,l$ respectively), with the union of these paths
containing a circuit surrounding $\mA$.  Here $\gamma_{g \lra h}$ denotes a path
from $g$ to $h$.  (Note that for general SLC $\mW$ and $\mS_k$, 
$\mR \bs (\mS_k \cup \mW)$ could include one or more
components consisting of a single bond with one endpoint in 
$V(\mW)$ and the other in $V(\mS_k)$, not surrounded by such lattice paths, but
that is not possible for in the present situation, due to our observation about
the outer surface inside $\Gamma_{\mR}$ for a ball.)
For lattice paths $\gamma,\tilde{\gamma}$
in $\mR$ having the same endpoints, we say that $\tilde{\gamma}$ is
\emph{directly obtainable from} $\gamma$ \emph{by contraction} if we can change
$\gamma$ to
$\tilde{\gamma}$ by one of the following two procedures:  (1) select a lattice
square
$Q$ such that two sides of $Q$ are consecutive bonds of $\gamma$, and replace
these two bonds with the other two sides of $Q$ (equivalently, replace the common
site of the two bonds with the opposite corner of $Q$, if we view the path as a
sequence of sites), or (2) select a lattice square $Q$ such that three sides of
$Q$ are consecutive bonds of $\gamma$, and replace them with the other side of $Q$
(i.e., shortcut the trip around $Q$.)  Again for paths $\gamma$ and
$\tilde{\gamma}$ having the same endpoints $g,h$, we say that $\tilde{\gamma}$ is
\emph{obtainable from} $\gamma$ \emph{by contraction in}
$\mR$ if there is a sequence of 
lattice paths in $\mR$, each from $g$ to $h$, starting with $\gamma$ and ending
with $\tilde{\gamma}$, each directly obtainable from the previous one by
contraction.  We may assume the paths $\gamma_{x \lra u}, \gamma_{x \lra v},
\gamma_{y \lra u}, \gamma_{y \lra v}$ are disjoint; if not, we replace the
starting site $x$ with the last common site of $\gamma_{x \lra u}$ and $\gamma_{x
\lra v}$ on the way to $u$ and $v$, respectively, from $x$, and similarly for $y$. 
Consider the two paths $\alpha_0 = \gamma_{x \lra u} \cup
\gamma_{y \lra u}$ and $\tilde{\alpha}_0 = \gamma_{x \lra v} \cup \gamma_{y \lra
v}$, each from
$x$ to $y$.  It is easy to see that there is a path $\alpha$ between these two
paths which is obtainable from both paths by contraction in $\mR$.  The various
paths, call them $\alpha_i$ and $\tilde{\alpha}_i$, obtained along the way 
from $\alpha_0$ and $\tilde{\alpha}_0$ to
$\alpha$ are each no
longer than the original paths 
$\alpha_0$ and $\tilde{\alpha}_0$, and every bond
between $\alpha_0$ and $\tilde{\alpha}_0$ is on one of the paths $\alpha_i$ or
$\tilde{\alpha_i}$.  The contraction aspect means that $\alpha_i,\tilde{\alpha_i}
\subset B_{\mR}(x,m) \cup \mW$, so we conclude that
$\mA \subset B_{\mR}(x,m) \cup \mW$, in contradiction to the definition of $\mA$. 
This establishes our claim that each
nontrivial component of $\mR \bs (\mS_k \cup \mW)$, and therefore also each
component of $\mR \bs (\mS_k
\cup \mW)$, contains exactly one segment of
$\Gamma_{\mR}$.  We thus have the following picture:  $\mR$ and $\mW$ are SLC, so
each component $\mC$ of $\mR \bs \mW$ includes exactly one segment, call it
$\gamma_{\mC}$, of 
$\Gamma_{\mR}$.  When we remove the set $\mS_k$ to obtain $\mR \bs (\mS_k \cup
\mW)$, the remaining portion $\mC \bs \mS_k$ of $\mC$ is all connected to
$\gamma_{\mC}$ and is the only component of $\mR \bs (\mS_k \cup \mW)$ which abuts
$\mS_k \cap \mC$.  This proves (c), and we have (a), (b), (c) while filling
$\Int(\Gamma_{\mR})$.

Since $\mV$ and $\mR \bs \mV$ are connected, $\mV \cap \Gamma_{\mR}$ is a single
segment, so after $\mV \cap \Int(\mR)$ is filled, we can fill $\mV \cap
\Gamma_{\mR}$ by starting at one end of this segment and proceeding to the other
end.  This keeps the part of
$\Gamma_{\mR}$ outside $\mS_k$ connected, and the above proof of (a), (b), (c)
remains valid (though now, in (c$^{\prime\prime}$), there may be more than two
components.)  As in Example \ref{Ex:slcprops}, this shows that
$(\mkR,\mkV,\mkS,\mkM,\mkN)$ is filling-compatible.
\end{example}

As Examples \ref{Ex:slcprops} and \ref{Ex:circbdd} show, it is sometimes necessary
to augment a natural collection, like the rectangles or the circuit-bounded sets,
to obtain classes
$\mkS, \mkV, \mkN$ having filling-compatibility, but this seems to be only a minor
technical obstactle to the applicability of our results to natural cases.

When dealing with the Ising model we will not have to restrict boundary
conditions, so in place of Definition \ref{D:enlarged} we can
use the following simpler ideas.  We say that a collection
$\mkL$ of finite subsets of $\ZZd$ has the \emph{approximate neighborhood property}
if for some $c > 1$, for every $\Lambda \in \mkL$ and $x \in \Lambda$, $\mkL$
includes a $c$-approximate $r$-neighborhood of $x$ in $\Lambda$   (A similar
property termed
\emph{inheriting} was used in
\cite{Al01fin}; the two are interchangeable for our purposes.)  We say that $\mkL$
is \emph{fillable} if for every
$\Lambda \in \mkL$ there exists and ordering $x_1,..,x_n$ of $\Lambda$ such that
for all $1 \leq k \leq n, \{x_1,..,x_k\} \in \mkL$.  Note that, in
contrast to the analogous properties (Definition \ref{D:enlarged}) for bond models,
here we need not incorporate collections of measures into the definitions, because
our collection of measures will always consist of all Ising models on all $\Lambda
\in \mkL$ with arbitrary boundary conditions.

We turn now to some definitions related to further conditions we will impose on the
collection $\mkM$ of measures.  For ordered
$J$, a probability measure $P$ on $J^{\Delta}$ for some finite
$\Delta$ is said to have
the \emph{FKG property} if $A, B$ increasing implies $P(A \cap B) \geq
P(A)P(B)$.  $P$ is said to satisfy the \emph{FKG lattice condition} if 
\begin{equation} \label{E:FKGlatt}
  P(\omega \vee \omega^{\prime})P(\omega \wedge \omega^{\prime})
    \geq P(\omega)P(\omega^{\prime}) \quad 
    \text{for all } \omega, \omega^{\prime},
\end{equation}
where $\vee$ and $\wedge$ denote the coordinatewise maximum and minimum,
respectively.  This implies that $P(\omega_{\Lambda} \in \cdot \mid 
\omega_{\Delta \backslash \Lambda} = \eta_{\Delta \backslash \Lambda})$ has the
FKG property for all $\Lambda$ and $\eta$ \cite{Ho}.  In a mild abuse of notation,
for a bond $e$ we write $\omega \vee e$ for the configuration taking value 1 at $e$
and agreeing with $\omega$ at all other bonds.  Then (\ref{E:FKGlatt}) is
equivalent to 
\[
  \frac{P(\omega \vee e)}{P(\omega)}\text{ is an increasing function of } \omega
    \text{ for each fixed } e.
\]
For $Q$ another probability
measure, we say that $P$ \emph{FKG-dominates}
$Q$ if for every nondecreasing function $f$ on
$J^{\Delta}$,
\[
  \int f dP \geq \int f dQ.
\]

We say that the collection $\mkM = \cup_{\mR \in \mkR}\ \mkM_{\mR}$ has 
\emph{uniform exponential decay of connectivity} if there exist
$C, \lambda > 0$ such that for every $\mR \in \mkR$ and
$P \in \mkM_{\mR}$,
\[
  P(x \lra y \text{ via a path in } \mR)  
    \leq C e^{-\lambda d_{\mR}(x,y)} \quad \text{for all } x, y \in V(\mR).
\]

The finite-volume analog of ratio weak mixing can be formulated for a
collection $\mkM = \cup_{\mR \in \mkR}\ \mkM_{\mR}$ of measures, as follows.  We
say that $\mkM$ has the \emph{ratio strong mixing property} if there exist $C,
\lambda > 0$ such that for all $\mR \in \mkR, P \in \mkM_{\mR}$ and $\mE, \mF
\subset \mR$, 
\begin{align} \label{E:rstrongmix}
  \sup &\left\{ \left| \frac{P(A \cap B)}{P(A)P(B)} - 1 \right| : A \in 
    \mathcal{G}_{\mathcal{E}}, B \in
    \mathcal{G}_{\mathcal{F}}, P(A)P(B) > 0 \right\} \\
  &\leq C \sum_{x \in V(\mathcal{E}),y \in V(\mathcal{F})}
    e^{-\lambda d_{\mR}(x,y)}, \notag
\end{align}
whenever the right side of (\ref{E:rstrongmix}) is less than 1.  This definition
was given in \cite{Al01fin} for the special case of a fixed bond percolation model
with some class of site or bond boundary conditions. 

A \emph{coupling} of two probability measures $P_{1}$ and $P_{2}$ on some set
$J^{\Delta}$ is a probability measure $\mathbb{P}$ on $J^{\Delta} \times
J^{\Delta}$ with marginals $P_{1}$ and $P_{2}$ (in order).  

For
$\Delta \subset \Lambda \subset \ZZd$ and $r > 0$ let
\[
  \Delta^r(\Lambda) = \{ x \in \Lambda: d_{\mB(\Lambda)}(x,\Delta) \leq r \},
    \qquad \Delta^r = \Delta^r(\ZZd).
\]
Similarly for $\mD \subset \mR \subset \mB(\ZZd)$ let
\[
  \mD^r(\mR) = \{ b \in \mR: d_{\mR}(b,\mD) \leq r \}, \qquad \mD^r =
    \mD^r(\mB(\ZZd)).
\]  

\section{Specific Models}

The \emph{FK model} (\cite{Fo1}, \cite{Fo2}, \cite{FK}; see also \cite{ACCN},
\cite{Gr95}) is a
graphical representation of the Potts model.  
For a configuration $\omega$ on $\mR \subset \mB(\ZZd)$, let
$K(\omega)$ denote the number of open clusters in $\omega$
which do not abut $\mR^c$.  For $p \in [0,1]$ and $q > 0$, the
FK model 
$P_{\mR,w}^{p,q}$ (without external fields) on $\mR$
with parameters $(p,q)$ and wired boundary condition is defined by the
weights
\begin{equation} \label{E:FKweight}
  W(\omega) = p^{|\omega|}(1 - p)^{|\mR| - |\omega|}q^{K(\omega)}.
\end{equation}
Here $|\omega|$ means the number of open bonds in $\omega$.  
Let $K(\omega \mid \rho)$ be the number
of open clusters of $(\omega \rho)$ which abut or intersect $\mR$.  The FK
model 
$P_{\mR,\rho}^{p,q}$ with
bond boundary condition $\rho$ is given by the weights in (\ref{E:FKweight}) with
$K(\omega)$ replaced by $K(\omega \mid \rho)$.  When $\rho = \rho^0$ or $\rho^1$ we
replace
$\rho$ with $f$ or $w$ in our notation.  The infinite-volume measures
\[
  P^{p,q}_{} = \lim_{\Lambda \nearrow \ZZd} P^{p,q}_{\Lambda,*}
\]
on $\{0,1\}^{\mB(\ZZd)}$ exist for $* = w$ or $f$ and are translation-invariant.
For $p$ below the percolation critical point $p_c(q,d)$ we have
$P^{p,q}_{w} = P^{p,q}_{f}$ so we omit the subscript.  For a summary of basic
properties of the FK model, see \cite{Gr95}.  In particular, for $q \geq 1$ the
FK model satisfies the FKG lattice condition, and we consider only these
values of $q$.

We also need to consider site boundary conditions, when we use the FK model as a
graphical representation of the Ising model.  Given $\Lambda \subset \ZZd$ and
$\eta \in \{-1,1\}^{\partial \Lambda}$ define
\[
  \partial \Lambda = \{ x \in \Lambda^c: x \text{ adjacent to } y \text{ for some }
    y \in \Lambda \},
\]
\begin{equation} 
  U(\Lambda,\eta) = \bigl\{\omega \in \{0,1\}^{\omB(\Lambda)}: \eta_{x}
    = \eta_{y} \text{ for every } x,y \in \partial \Lambda \text{ for which }
    x \lra y \text{ in } \omega \bigr\}. \notag
\end{equation}
The FK model $P^{p,q}_{\omB(\Lambda),\eta}$ with
site boundary condition $\eta$ is given by the weights in (\ref{E:FKweight}),
multiplied by $\delta_{U(\Lambda,\eta)}(\omega)$.  

For the FK model with external fields
$h_{i}$, $i = 1,\ldots,q$ and free boundary, 
the factor $q^{K(\omega)}$ in the weight $W(\omega)$ is
replaced by 
\begin{equation} \label{E:FKexternal}
  \prod_{C \in \mK(\omega)} \left((1-p)^{h_{1}s(C)} + (1-p)^{h_{2}s(C)} +
  \ldots + (1-p)^{h_{q}s(C)}\right),
\end{equation}
where $\mK(\omega)$ is the set of open clusters in $\mR$ in the
configuration $\omega$ and
$s(C)$ denotes the number of sites in the cluster $C$.  The parameters are then
$(p,q,\{h_{i}\})$; $q$ must be an integer, and
we may omit $\{h_{i}\}$ when all external fields are 0.  The percolation critical
point is denoted $p_c(q,d,\{h_i\})$.  We need only consider 
$0 = h_{1} \geq h_{2} \geq \ldots \geq h_{q}$, so we henceforth assume this in
our notation.  Species $i$ is called \emph{stable} if
$h_{i}$ is maximal, i.e. $h_{i} = h_{1} = 0$.  For bond
boundary conditions $\rho$ we replace (\ref{E:FKexternal}) with
\begin{equation} \label{E:FKextbond}
  \prod_{C \in \mK(\omega \mid \rho)} \left((1-p)^{h_{1}s(C)} +
  (1-p)^{h_{2}s(C)} +
  \ldots + (1-p)^{h_{q}s(C)}\right),
\end{equation}
where $\mK(\omega \mid \rho)$ is the set of finite open 
clusters of $(\omega\rho)$ which
intersect $V(\mB)$.  For
general site boundary conditions $\eta$ for the model on $\omB(\Lambda)$ the factor
(\ref{E:FKexternal}) is multiplied by 
\begin{align} \label{E:FKextsite}
  \prod_{C \in \mK_{int}(\omega)} &\left((1-p)^{h_{1}s(C)} +
    (1-p)^{h_{2}s(C)} +
    \ldots + (1-p)^{h_{q}s(C)}\right) \\
  &\times \quad
    \prod_{C \in \mK_{\partial}(\omega)} (1-p)^{h_{i(C)}s(C)} \quad \times \quad
    \delta_{U(\Lambda,\eta)}(\omega), \notag
\end{align}
where $\mK_{\partial}(\omega)$ (respectively $\mK_{int}(\omega)$) is the set of
clusters in the configuration
$\omega$ which do (respectively don't) intersect 
$\partial \Lambda$ and $i(C)$ is
the species for which $\eta_{x} = i$ for all $i \in \partial \Lambda \cap C$. 
(The existence of such an $i$ is forced by the event $U(\Lambda,\eta)$.)  

If $\omB(\Lambda)^c$ is
connected then for stable $i$ the wired boundary
condition is equivalent to the all-$i$ site boundary condition.   
For $q \geq 1$, the FK model with external fields satisfies the FKG lattice
condition, under any bond boundary
condition.  We say that $\rho \in
\{0,1\}^{\mR^c}$ is a \emph{unique-cluster bond boundary condition} if
all open bonds in $\rho$ are part of one cluster.
We have seen (Figure \ref{F:markov}) that in the absense of external fields, an FK
model
$P_{\mR,\rho}^{p,q}$ need not in general have the Markov property for blocking
sets, but we have the following sufficient conditions.

\begin{lemma} \label{L:markov}
  Let $\mR \subset \mB(\ZZd)$ and consider an FK measure
  $P_{\mR,\rho}^{p,q,\{h_i\}}$.  

  (i) Suppose $\rho$ is a unique-cluster bond boundary
  condition, and suppose that either (a) there
  are no external fields, (b) there is a  unique nonsingleton cluster in $\rho$
  and this cluster is infinite, or (c) $\rho = \rho^0$. 
  Then $P_{\mR,\rho}^{p,q,\{h_i\}}$ has
  the Markov property for blocking sets.  

  (ii) If $\rho$ is arbitrary, there are no external fields,
  $\mE \subset \mR$ and each component of $\mR \bs \mE$ abuts at most
  one nonsingleton cluster of $\rho$, then $P_{\mR,\rho}^{p,q,\{h_i\}}$ has the 
  Markov property for sets blocking $\mE$.
\end{lemma}
\begin{proof}
We first prove (i).  The FK weight can be written as a product over
clusters,
\[
  \prod_{C \in \mK(\omega \mid \rho)} \left( \frac{p}{1-p} \right)^{b(C)}
  \left((1-p)^{h_{1}s(C)} + (1-p)^{h_{2}s(C)} +
  \ldots + (1-p)^{h_{q}s(C)}\right),
\]
where $s(C)$ and $b(C)$ are the number of sites and bonds, respectively, in the
cluster $C$.  Let $C_u$ denote the unique nonsingleton cluster in $\rho$, when
this exists.  Let $(\mX,\mY,\mZ)$ be a blocking partitiion of $\mR$, and suppose
$\omega = 0$ on $\mY$.  Note that a group of clusters $C_1,..,C_n$ of $\omega$ may
be part of the same cluster, say  $\hat{C}$, in $(\omega\rho)$, if $C_1,..,C_n$
are connected together by
$C_u$; in particular this can occur with some of the $C_i$'s on each side of the
blocking set $\mY$.  However, under (a), (b) or (c), the weight of
$\hat{C}$ factors into a product of a weight for each $C_i$, 
and thus the clusters on
each side of the blocking set occur independently, yielding the Markov property.
The proof of (ii) is similar.
\end{proof}

Lemma \ref{L:markov} is part of what requires us to use augmented collections
instead of using a single $\mkR$ and $\mkM$ throughout; in Example
\ref{Ex:slcprops}, for example, we restrict $\mkM$ to free and wired boundary
conditions to guarantee the Markov property for blocking sets, but such a
restriction would be unnecessary and technically awkward for $\mkM_w^+(\mkN,\mkR)$,
which does not need the Markov property in our proofs.

The following facts about the FK model are known for $d=2$. For $q = 1, q = 2,$ and
$q \geq 25.72$, we have $p_{c}(q,2) = 
\tfrac{\sqrt{q}}{1 + \sqrt{q}}$ \cite{LMR}, and the connectivity decays 
exponentially for all $p < p_{c}(q,2)$ \cite{Gr97}.  This is believed to be true
for all $q$; for $2 < q < 25.72$ the connectivity is known to
decay exponentially at least for all $p < \tfrac{\sqrt{q-1}}
{1 + \sqrt{q-1}}$, and analogous results hold for other planar lattices
\cite{AlARC}. For general $q \geq 1$,  if the
connectivity decays exponentially then the model has the ratio weak  mixing
property \cite{Al98mix}.  (This result is actually given 
assuming a nonnegative external field applied to
at most one species, but the proof carries over without change to arbitrary
external fields; the necessary FKG property is 
proved in \cite{BBCK}.) 

As shown in \cite{ES}, for $\beta$ given by $p = 1 -
e^{-\beta}$, a configuration of the Ising model on $\Lambda$ with
boundary condition $\eta$ at inverse temperature $\beta$ can be obtained from a
configuration
$\omega$ of the FK model at $(p,2)$ with site boundary condition $\eta$, by
choosing a label for each cluster of
$\omega$ independently and uniformly from $\{-1,1\}$; this
\emph{cluster-labeling construction} yields a joint site-bond configuration for
which the sites are an Ising model and the bonds are an FK model.  When the
parameters are related in this way, we call the Ising and FK models
\emph{corresponding}.  Alternatively, if one selects an Ising configuration
$\sigma_{\Lambda}$ and does independent percolation at density $p$ on the set of
bonds
\[
  \{\langle xy \rangle \in \omB(\Lambda): (\sigma\eta)_{x} = 
    (\sigma\eta)_{y} \},
\]
the resulting bond configuration is a realization of the corresponding FK model. 
We call this the \emph{percolation construction} of the FK model.  

For the Ising model at inverse temperature $\beta < \beta_{c}(d)$, for $p = 1 -
e^{-\beta}$ and for the FK model without external fields at
$(p,2)$, the covariance in the Ising model and the connectivity in the FK model
are related by
\begin{equation} \label{E:covarconn}
  \cov(\sigma_0,\sigma_x) = P(0 \lra x);
\end{equation}
see \cite{ACCN} or \cite{Gr95}.  Thus exponential decay of connectivities 
in the FK model is equivalent to exponential decay of correlations in
the corresponding Ising model.  Further, 
$\beta_{c}(d)$ and the percolation critical point
$p_{c}(2,d)$ of the FK model are related by
\[
  p_{c}(2,d) = 1 - e^{-\beta_{c}(d)};
\]
again see \cite{ACCN} or \cite{Gr95}.  For $h \neq 0$ we make this a definition,
that is, $\beta_c(d,h)$ is defined by
\[
  p_{c}(2,d,h) = 1 - e^{-\beta_{c}(d,h)},
\]
where $p_c(2,d,h)$ is the percolation threshhold of the corresponding FK
model.  (The notation is not meant to imply that $\beta_c(d,h)$ is a true critical
point.)

\section{Statements of Main Theorems}
All proofs appear in Section 4.
Our first theorem covers bond percolation models in finite volumes $\mR$.  Note
that as discussed in the introduction, one of the two events $A,B$ is restricted to
occur somewhere on a fixed set
$\mD \subset \mR$, and the required separation $r$ depends on the size of $\mD$. 
The location of the other event is unrestricted and in particular this location may
also be a part of $\mD$.  

\begin{theorem} \label{T:finitevol}
  Let $\mkR, \mkV, \mkS, \mkM$ be as in the standard bond percolation setup and
  let $\mkN$ be a neighborhood collection for $\mkS$.
  Suppose that
  $(\mkR, \mkV, \mkS, \mkM,\mkN)$ is filling-compatible, each 
  measure in $\mkM$ satisfies the FKG lattice condition,
  each measure in $\mkM \cup \mkM^+(\mkS,\mkR)$
  has the Markov property for blocking sets, and 
  $\mkM_w^+(\mkN,\mkR)$ has uniform exponential decay of connectivity.
  There exist $c_{i},\ep_i$ 
  such that for all $\mR \in \mkR, P \in \mkM_{\mR}, \mV \in \mkV_{\mR}$, and
  $r \geq c_{7} \log |\mV|$, all $\mD$ with $\mD^r(\mR) \subset \mV$
  and all increasing or decreasing events $A
  \in \mG_{\mR}$ 
  and $B \in \mathcal{G}_{\mathcal{D}}$,
  \[
    P(A \circ_{r} B) \leq (1 + c_{8}e^{-\ep_{4}r}) P(A)P(B).
  \]
\end{theorem}

If in Theorem \ref{T:finitevol} we only assume filling-compatibility at a
particular scale $R>1$, then the conclusion is valid provided we further restrict
to $r \leq R$.

In the case of the FK model on $\mB(\ZZd)$, Theorem \ref{T:finitevol} will yield
the next theorem.  Site boundary conditions cannot be allowed in Theorem
\ref{T:FKZd}, as the corresponding measures need not in general have the FKG
property.  Multiple-cluster bond boundary conditions cannot be allowed due to the
phenomenon, dubbed \emph{tunneling} in
\cite{Al01fin}, that such boundary conditions may create long-range dependencies,
even when the locations of the events $A,B$ are nonrandom.  In other words, the
strong mixing property may fail.  In fact we restrict ourselves to wired boundary
conditions, in order to obtain the filling-compatible property in a
straightforward way, but one can presumably allow more general unique-cluster bond
boundary conditions.

We say that an event $A \in \mG_{\omB(\ZZd)}$ is \emph{locally-occurring} if
$\omega
\in A$ implies that $A$ occurs in $\omega$ on some finite set of bonds; for
site models an analogous definition is made for $A \in \mH_{\ZZd}$.

\begin{theorem} \label{T:FKZd}
  Let $\mkR$ be the collection of all lattice rectangles, and $\omkR$
  the class of all approximate lattice rectangles, in $\mB(\ZZd)$,
  and for $\mR \in \mkR$ let $\mkV_{\mR}$ be the collection of all lattice
  rectangles $\mV \subset \mR$ with $\mR \bs \mV$ connected. 
  Let $P = P^{p,q,\{h_i\}}$ be an FK model on $\mB(\ZZd)$.  Suppose $P$ has
  uniform exponential decay of connectivity for the class $\omkR$ with wired
  boundary conditions.  There exist $c_i,\ep_i$ such that the following hold.

  (i) For all $\mR \in \mkR$ and $\mV \in \mkV_{\mR}$, for all $\mD$ with
  $\mD^r(\mR) \subset \mV$,
  for all increasing or decreasing events $A \in \mG_{\mR}, B \in \mG_{\mD}$, 
  for all $r \geq c_{9} \log |\mV|$, and for $*=w$ or $f$,
  \[
    P_{\mD,*}(A \circ_{r} B) \leq 
      (1 + c_{10}e^{-\ep_{5}r}) P_{\mD,*}(A)
      P_{\mD,*}(B). 
  \]

  (ii) For all increasing or decreasing events $A,B$ with 
  $A$ locally-occurring and $B \in \mG_{\mD}$, and
  for all $r \geq c_{9} \log \diam(\mD)$,
  \[
    P(A \circ_{r} B) \leq (1 + c_{10}e^{-\ep_{5}r}) P(A)P(B). 
  \]
\end{theorem}

A sufficient condition for the uniform exponential decay hypothesized in Theorem
\ref{T:FKZd} is that $p$ be below the percolation critical point $p_c(1,d)$ for
independent percolation on $\mB(\ZZd)$.  The FK model at $(p,q,\{h_i\})$ is
FKG-dominated by independent percolation at density $p$, since $q \geq 1$, and
independent percolation at every density $p < p_c(1,d)$ has (uniform) exponential
decay of connectivity \cite{Me}.  

Using Lemma \ref{L:markov}, one could presumably extend Theorem \ref{T:FKZd} beyond
free and wired boundary conditions to general unique-cluster bond boundary
conditions, by using a more elaborate filling algorithm than the one in Example
\ref{Ex:slcprops}.  If $\mC$ is the unique nonsingleton cluster in a boundary
condition $\rho$, one would have to keep $\mC \cup (\mR \bs \mS)$ connected as the
filling proceeds, so that the effective boundary condition on $\mS$ when we
condition $P_{\mR,\rho}$ on the event $\omega_{\mR \bs \mS} = \rho_{\mR \bs
\mS}^1$ still has a unique nonsingleton cluster.  But since we have no specific
example as motivation for undertaking the additional technicalities, we will not
do so here.  

One can readily extend Theorem \ref{T:FKZd}(ii) to allow $A$ to be a limit in an
appropriate sense of locally-occurring events, but considering completely general
increasing $A$ creates technical difficulties; lacking again a motivating example
we have not attempted to surmount these.

\begin{remark} \label{R:irregR}
In Theorem \ref{T:FKZd}, as will be apparent from the proof, we need not require
$\mR$ to be a lattice rectangle if the following condition is satisfied: 
$\mR$ is connected and there exists a lattice rectangle $\mV^{\prime\prime}$ such
that $\mD^r \subset
\mV^{\prime\prime}
\subset \mR$.  We then let $\mV, \mV^{\prime}$ be lattice
rectangles which have the same center as $\mV^{\prime\prime}$ but are respectively
$r$ and $r/2$ units shorter in each direction, meaning $\mV \supset \mD^{r/2}$;
we use
$r/4$ in the proof in place of $r$.  When $\mV$ is being filled, 
the relevant sets $\mS_k \cup \mW$ as in Definition
\ref{D:enlarged} are contained in $\mV^{\prime}$.
When we condition in a way that forces 
$\omega_{\mR \bs \mV^{\prime}} = \rho_{\mR \bs \mV^{\prime}}^i$ for $i = 0$ or 1,
the effective boundary condition for configurations on $\mV^{\prime}$ is free
or wired, respectively, regardless of what is outside $\mV^{\prime\prime}$.  
\end{remark}

For $d=2$ we have the following stronger result; we need not explicitly
assume
\emph{uniform} exponential decay of connectivity because, in SLC sets, it follows
from the usual infinite-volume exponential decay of connectivity \cite{Al01fin}.

\begin{theorem} \label{T:FKZ2}
  Let $\mkR$ be the collection of all circuit-bounded subsets of
  $\mB(\ZZ)$, let $\mkV_{\mR}$ be as in Example \ref{Ex:circbdd} for $\mR \in
  \mkR$, and let $P = P^{p,q,\{h_i\}}$ be an FK model on $\mB(\ZZ)$.  Suppose
  $p < p_c(q,d,\{h_i\})$ and $P$ has
  exponential decay of connectivity (in infinite volume.) 
  There exist $c_j$ such that for all $\mR \in
  \mkR$, all $\mV \in \mkV_{\mR}$, all $r \geq c_{11} \log |\mV|$,
  and all $\mD \subset \mR$ with $\mD^r(\mR) \subset \mV$,
  for all increasing or decreasing events $A \in \mG_{\mR}, B \in \mG_{\mD}$, 
  and for $*=w$ or $f$, 
  \begin{equation} \label{E:sepocc}
    P_{\mR,*}(A \circ_{r} B) \leq (1 + c_{12}e^{-\ep_{6}r})
      P_{\mR,*}(A) P_{\mR,*}(B).
  \end{equation}
\end{theorem}

\begin{remark} \label{R:SLCC}
Theorem \ref{T:FKZ2} extends straighforwardly to the case in which $\mR$ has
multiple components, each circuit-bounded, using the fact that under free and wired
boundary conditions, the configurations on the various components are independent.
\end{remark}

It is possible to allow general site and bond boundary conditions for the FK model
if we resrict one of the two events to occur well-separated from the
boundary, specifically the event restricted to occur on a particular
$\mD$.  This means we assume $\mD^r \subset \mV$ instead of $\mD^r(\mR) \subset
\mV$.

\begin{theorem} \label{T:farbdry}
  Let $\mkR$ be the collection of all circuit-bounded subsets of
  $\mB(\ZZ)$, let $\mkV_{\mR}$ be as in Example \ref{Ex:circbdd} for $\mR \in
  \mkR$, and let $P = P^{p,q,\{h_i\}}$ be an FK model on $\mB(\ZZ)$.  Suppose
  $p < p_c(q,d,\{h_i\})$ and $P$ has
  exponential decay of connectivity (in infinite volume.) 
  There exist $c_j$ such that for all $\mR \in
  \mkR$, all $\mV \in \mkV_{\mR}$, all $r \geq c_{11} \log |\mV|$,
  and all $\mD \subset \mR$ with $\mD^r \subset \mV$,
  for all increasing or decreasing events $A \in \mG_{\mR}, B \in \mG_{\mD}$,  
  and for all site or bond boundary conditions $\rho$,
  \[
    P_{\mR,\rho}(A \circ_{r} B) \leq (1 + c_{14}e^{-\ep_{7}r})
      P_{\mR,\rho}(A) P_{\mR,\rho}(B). 
  \]
\end{theorem}

The last two of our main theorems cover the Ising model.  An \emph{absorbing
sequence} in $\ZZd$ is an increasing sequence of subsets whose union is $\ZZd$.

\begin{theorem} \label{T:Ising}
  Let $\mu = \mu^{\beta,h}$ be an Ising model on $\ZZd$, let $\mkL$ be a
  collection of finite subsets of $\ZZd$ which is fillable and has the neighborhood
  component property.  Suppose the corresponding FK model has uniform exponential
  decay of connectivity for the class $\{ \mB(\Lambda): \Lambda \in \mkL \}$ 
  with wired boundary conditions.
  There exist $c_i,\ep_i$ such that for all finite $\Delta \subset
  \ZZd$, the following hold.

  (i) For all $\Theta,\Lambda \in \mkL$ with $\Theta \subset
  \Lambda$ and $|\Theta| \geq 2$, for all $r \geq c_{15} \log |\Theta|$,
  for all boundary conditions $\eta$, for all $\Delta \subset \ZZd$ with
  $\Delta^r(\Lambda) \subset \Theta$, and
  for all increasing or decreasing events $A \in \mH_{\Lambda}, B \in
  \mH_{\Delta}$, 
  \[
    \mu_{\Lambda,\eta}(A \circ_{r} B) \leq (1 + c_{16}e^{-\ep_{8}r}) 
      \mu_{\Lambda,\eta}(A)\mu_{\Lambda,\eta}(B). 
  \]

  (ii) Assume $\mkL$ contains an absorbing sequence.  Then for all $\Theta \in
  \mkL$ and all $r \geq c_{15} \log |\Theta|$, for all $\Delta \subset \ZZd$ with
  $\Delta^r \subset \Theta$, and
  for all increasing or decreasing events $A,B$ with $B \in \mH_{\Delta}$
  and $A$ locallly-occurring,  
  \[
    \mu(A \circ_{r} B) \leq (1 + c_{16}e^{-\ep_{8}r}) \mu(A)\mu(B). 
  \]
\end{theorem}

We now specialize to SLC subsets in two dimensions.  As noted above, it is proved
in
\cite{Al01fin} that the hypothesis in Theorem \ref{T:Ising2} of uniform exponential
decay of connectivity in the corresponding FK model is satisfied whenever that FK
model has exponential decay of connectivity in infinite volume.  If $h=0$, then by
(\ref{E:covarconn}) this exponential decay of connectivity
in infinite volume holds whenever there is a unique Gibbs distribution and this
distribution has exponential decay of correlations, i.e. whenever
$\beta < \beta_c(2,0)$ \cite{ABF}.

\begin{theorem} \label{T:Ising2}
  Let $\mu = \mu^{\beta,h}$ be an Ising model on $\ZZ$ and let $\mkL$ be the
  classs of all finite SLC subsets of $\ZZ$ with arbitrary boundary condition. 
  Suppose that either (a) $\beta <
  \beta_c(2,0)$ and $h = 0$, (b) $\beta > \beta_c(2,0)$ and $h \neq 0$, (c) $\beta
  < \beta_c(2,h)$ and the corresponding
  FK model has exponential decay of connectivities (in infinite volume), or (d)
  $\beta > \beta_c(2,h)$ and the corresponding FK model has exponential decay of
  dual connectivities (in infinite volume).  
  There exist $c_i,\ep_i$ such that for all $\Delta \subset
  \ZZd$ with $|\Delta| \geq 3$, the following hold.

  (i) For all $\Lambda \in \mkL$ with $\Delta \subset \Lambda$,
  for all boundary conditions $\eta$, 
  for all increasing or decreasing events $A \in \mH_{\Lambda}, B \in
  \mH_{\Delta}$, and
  for all $r \geq c_{17} \log \diam_{\mB(\Lambda)}(\Delta)$,
  \[
    \mu_{\Lambda,\eta}(A \circ_{r} B) \leq (1 + c_{18}e^{-\ep_{9}r}) 
      \mu_{\Lambda,\eta}(A)\mu_{\Lambda,\eta}(B). 
  \]

  (ii) For all increasing or decreasing events $A,B$ with $B \in \mH_{\Delta}$ and
  $A$ locally-occurring, and
  for all $r \geq c_{17} \log \diam(\Delta)$,
  \[
    \mu(A \circ_{r} B) \leq (1 + c_{18}e^{-\ep_{9}r}) \mu(A)\mu(B). 
  \]
\end{theorem}

\section{Proofs}
We begin with the proof of Theorem \ref{T:finitevol}.  
Since $|\mV| \geq 2$, we need only consider ``sufficiently large''$r$; we do this
tacitly throughout.
Our proof will be based on an elaboration of the
``bond-splitting'' proof, given by van den Berg and
Fiebig in \cite{vdBF}, of the van den Berg-Kesten
inequality
\cite{vdBK}, so we begin with a brief review of the basic idea of \cite{vdBF}.  For
independent percolation on a finite set
$\mR$ of bonds, one may take $\mathcal{V} \subset \mR$ and
``split'' each bond in
$\mathcal{V}$ into, say, a left and a right bond; the left and right bonds receive
open/closed states independently.  For increasing or decreasing events $A$ and
$B$, one considers the event that $A$ and $B$ occur disjointly, with $A$ occuring
in the configuration of unsplit and left bonds, and $B$ occurring in the
configuration of unsplit and right bonds.  One shows that splitting an additional
bond never decreases the probability of this form of disjoint occurrence.  When
all bonds are split, $A$ and $B$ become independent, yielding the inequality.

For dependent models this does not work in general.  For 
example, if one considers the FK model on a graph with some set
$\mathcal{S} \subset \mV$ of split bonds, the marginal distribution
of the configuration on the unsplit and left bonds is not the same as the
distribution of the original model on the fully-unsplit graph.  Instead, for a
bond percolation model $P$ on a set $\mR$ of bonds,
writing $\mathcal{T}$ for $\mathcal{R}
\backslash \mathcal{S}$, we consider $\mathcal{S}$-\emph{split
configurations} $(\omT,\omS,\omSt) \in
\{0,1\}^{\mathcal{T}} \times \{0,1\}^{\mathcal{S}} \times \{0,1\}^{\mathcal{S}}$
under the probability measure 
\[
  P_{\mathcal{S}}(\omT,\omS,\omSt) = P(\omT)P(\omS
  \mid \omT)P(\omSt \mid \omT),
\]
which we call the $\mathcal{S}$-\emph{split measure}.  Note $\omS,\omSt$ are
conditionally independent given $\omT$.  For $A, B \subset
\{0,1\}^{\mR}$ and $(\omT,\omS,\omSt) \in
\{0,1\}^{\mathcal{T}} \times \{0,1\}^{\mathcal{S}} \times \{0,1\}^{\mathcal{S}}$,
we say that $A$ and $B$ \emph{occur}
$\mathcal{S}$-\emph{split at separation} $r$ in $(\omT,\omS,\omSt)$ if there
exist $\mathcal{E},\mathcal{F} \subset \mR$ with
$d_{\mR}(\mathcal{E},\mathcal{F}) \geq r$ such that $A$ occurs on $\mathcal{E}$ in
$(\omT,\omS)$ and $B$ occurs on $\mathcal{F}$ in $(\omT,\omSt)$.  We denote this
event by $A \circ_{r,\mathcal{S}} B$.  Let $\mV = \{b_1,..,b_n\}$ be a
filling sequence.  If we can show that for some $c_{19},\ep_{10}$, for each $k
\leq n$, for $\mS = \{ b_1,..,b_{k-1} \}$ and $e = b_k$, we have
\begin{equation} \label{E:inducstep}
  P_{\mathcal{S}}(A \circ_{r,\mathcal{S}} B) \leq (1 + c_{19} e^{-\ep_{10} r})
    P_{\mathcal{S} \cup \{e\}}(A \circ_{r,\mathcal{S} \cup \{e\}} B),
\end{equation}
then we obtain
by iterating (\ref{E:inducstep}) that
\begin{align} \label{E:induc}
  P(A \circ_{r} B) &= P_{\phi}(A \circ_{r,\phi} B) \\
  &\leq  (1 + c_{19} e^{-\ep_{10} r})^{|\mV)|}
    P_{\mV}(A \circ_{r,\mV} B). \notag 
\end{align}
Now since $B$ occurs only on $\mD$, and $\mD \subset \mV$,
\begin{align} \label{E:allsplit}
  P_{\mV}&(A \circ_{r,\mV} B) \\
  &\leq P_{\mV}(A \text{ occurs in }
    (\omega_{\mR \bs \mV},\omega_{\mV})) P_{\mV}(A \circ_{r,\mV}
    B \mid A \text{ occurs in }
    (\omega_{\mR \bs \mV},\omega_{\mV}))\notag \\
  &\leq P(A) \sup_{\rho_{\mR \bs \mV}}
    P(B \mid \omega_{\mR \bs \mV} = \rho_{\mR \bs \mV}). \notag 
\end{align}
Finally we can apply Proposition \ref{P:strongmix} below to the collection $\mkS_0
= \{ \{ \mR
\}: \mR \in \mkR \}$, noting that $\mkM^+(\mkS_0,\mkR) = \mkM$ and
$d_{\mR}(\mD,\mR \bs \mV) \geq r$, to conclude that ratio strong mixing applies,
and the right side of (\ref{E:allsplit}) is bounded by
\begin{align} \label{E:allsplit2}
  &P(A)P(B) \left( 1 + \sum_{x \in V(\mD),
    y \in V(\mR \bs \mV)} c_{20}e^{-\ep_{11}d_{\mR}(x,y)} \right) \notag \\
  &\quad \leq P(A)P(B) \left( 1 + c_{21}e^{-\ep_{12}r} \right). \notag
\end{align}
Since
\[
  (1 + c_{19} e^{-\ep_{10} r})^{|\mV|}
  \leq 1 + c_{22} e^{-\ep_{13} r}, 
\]
this will complete the proof.

The main difficulty in proving (\ref{E:inducstep}) is that the properties which
can be established for the model
$P$, particularly weak mixing, do not immediately carry over to the
$\mathcal{S}$-split measure.  

Our proof of (\ref{E:inducstep}) will involve a coupling construction which we now
describe.  Fix $A$ and $B$.  If one of $A,B$ is increasing and the other is
decreasing, then since $P$ has the FKG property, the separated-occurrence
inequality is trivial:
\[
  P(A \circ_{r} B) \leq P(A \cap B) \leq P(A)P(B),
\]
so we may assume $A,B$ are both increasing or both decreasing.
Let $\mathcal{U} = \mathcal{T} \backslash \{e\}$.
Suppose that for each $\zeU \in \{0,1\}^{\mathcal{U}}$ 
we have a measure $\hat{P}_{\zeU}$ on
$\{0,1\}^{\mathcal{S}} \times \{0,1\}^{\mathcal{S}}$, which is a coupling of
$P(\omS \in \cdot \mid \omU = \zeU, \ome = 1)$ and $P(\omS \in \cdot \mid \omU
= \zeU, \ome = 0)$ satisfying
\begin{equation} \label{E:FKGcoup}
  \hat{P}_{\zeU}\left(\{ (\omS^{1},\omS^{0}): \omS^{1} 
    \geq \omS^{0} \}\right) = 1.
\end{equation}
A coupling for which (\ref{E:FKGcoup}) holds is called an
\emph{FKG coupling}.  An FKG coupling always exists when the first measure
FKG-dominates the second, as is the case here; see \cite{Ho}. Define
$\mathbb{P}_{\mathcal{S},e}$
on $\{0,1\}^{\mathcal{U}} \times \{0,1\}^{2}
\times (\{0,1\}^{\mathcal{S}})^{4}$ by
\begin{align} \label{E:PSdef}
  \mathbb{P}_{\mathcal{S},e}&(\zeU,\zte,
    \zeet,\zeS^{1},\zeS^{0},\zeSt^{1},\zeSt^{0}) \\
  &= P(\omU = \zeU) P(\ome = \zte \mid \omU = \zeU)
    P(\ome = \zeet \mid \omU = \zeU) \hat{P}_{\zeU}\bigl((\zeS^{1},\zeS^{0})\bigr)
    \hat{P}_{\zeU}\bigl((\zeSt^{1},\zeSt^{0})\bigr). \notag
\end{align}
To explain, for events $A,B \in \{0,1\}^{\mR}$ the measure
$\mathbb{P}_{\mathcal{S},e}$ arises in the following construction.  
First choose $\omU$ under the measure $P$.  Then choose
what we will call the $A$ \emph{pair} $(\omS^{1},\omS^{0})$ 
using the coupling measure
$\hat{P}_{\omU}$, and then independently (given $\omU$) the
$B$ \emph{pair} $(\omSt^{1},\omSt^{0})$ 
again using the coupling measure
$\hat{P}_{\omU}$.  We refer to $\omS^{1}$ or $\omSt^{1}$ as the
\emph{top layer}, and to $\omS^{0}$ or $\omSt^{0}$ as the \emph{bottom layer}, in
its respective pair.  We then choose $\ome,\omet$ independently (given $\omU$)
under the measure $P(\ome = \cdot \mid \omU = \zeta_{\mU})$ and use these values to
determine which layer to use from each pair in forming an $(\mS \cup \{e\})$-split
configuration.  If, for example, $\ome = 1$ and
$\omet = 0$, we form an $(\mathcal{S} \cup \{e\})$-split configuration out 
of $\omU, \ome,
\omet$, the top layer of the $A$ pair and the bottom layer of the $B$ pair,
and we may look for the event $A \circ_{r,\mathcal{S} \cup \{e\} } B$ in this
configuration. Other values of $\ome,\omet$ give corresponding different choices
of top or bottom layers to use in the $(\mathcal{S} \cup \{e\})$-split
configuration.  Note that $\omS^{i}$ and
$\omSt^{j}$ are conditionally independent given $\omU$, for each $i,j = 0,1$;
 from this it is easy to see that the constructed configuration $\bigl(
\omU,(\ome,\omS^{\ome}),(\ome,\omSt^{\omet}) \bigr)$ has the $(\mathcal{S} \cup
\{e\})$-split measure as its distribution.  By contrast, as a different
construction, instead of using
$\ome,\omet$ to choose a layer in each pair, we may use a single variable, say
$\ome$, and use it to choose a layer in both pairs; that is, we use the top layer
in both pairs if $\ome = 1$, and the bottom layer in both pairs if $\ome = 0$.
The resulting configuration $\bigl( (\omU,\ome),\omS^{\ome},\omSt^{\ome} \bigr)$
has the $\mathcal{S}$-split measure as its distribution.  The split-occurrence
events corresponding to these two constructions are
\begin{align} 
  &C(A,B,r,\mathcal{S} \cup \{e\}) = \notag \\
  &\qquad \bigl\{
    (\omU,\ome,\omet,\omS^{1},\omS^{0},\omSt^{1},
    \omSt^{0}):  \bigl( \omU,(\ome,\omS^{\ome}),(\omet,\omSt^{\omet}) \bigr) 
    \in A \circ_{r,\mathcal{S} \cup \{e\}} B \bigr\}, \notag \\
  &C(A,B,r,\mathcal{S}) = \notag \\
  &\qquad \bigl\{
    (\omU,\ome,\omet,\omS^{1},\omS^{0},\omSt^{1},
    \omSt^{0}):  \bigl( (\omU,\ome),\omS^{\ome},\omSt^{\ome} \bigr) 
    \in A \circ_{r,\mathcal{S}} B \bigr\}. \notag 
\end{align}
Thus we have
\begin{align} \label{E:Cprobs}
  \mathbb{P}_{\mathcal{S},e}(C(A,B,r,\mathcal{S} \cup \{e\})) &= 
    P_{\mathcal{S} \cup \{e\} }(A \circ_{r,\mathcal{S} \cup \{e\} } B), \\
  \mathbb{P}_{\mathcal{S},e}(C(A,B,r,\mathcal{S})) &= 
    P_{\mathcal{S}}(A \circ_{r,\mathcal{S}} B). \notag 
\end{align}

For fixed $\zeU,\zeS^{1},\zeS^{0},\zeSt^{1},\zeSt^{0}$ we may ask,
which of these constructions gives the greater probability for split
occurrence at separation $r$?  (The proof of the van den Berg-Kesten inequality
is based on the fact that in the independent context the first construction---with
two separate variables
$\ome,\omet$---always gives the higher probability; in that context
only one layer is needed for each of $A,B$ instead of two.)  To approach this
question in the present context, note that at least one of the sets
$\mathcal{E},\mathcal{F}$ where $A,B$ occur must be outside $\{e\}^{r/3}$.
Suppose now that $A,B$ are increasing; the decreasing case is
similar.  For $i, j = 0,1$ set
\[
  C_{ij} = \bigl\{
    (\omU,\omS^{1},\omS^{0},\omSt^{1},
    \omSt^{0}):  \bigl( \omU,(i,\omS^{i}),(j,\omSt^{j}) \bigr) 
    \in A \circ_{r,\mathcal{S} \cup \{e\}} B \bigr\}.
\]
This is the event that, loosely, ``if the $i$ and $j$ layers are chosen 
in the $A$ and $B$ pairs respectively, then $A
\circ_{r,\mathcal{S} \cup \{e\}} B$ will occur.''  We add a superscript $A$ or
$B$ to
$C_{ij}$ to designate which event occurs outside $\{e\}^{r/3}$, so that for
example
$C_{ij}^A$ is the event that ``if the $i$ and $j$ layers are chosen, then $A
\circ_{r,\mathcal{S} \cup \{e\}} B$ will occur with $A$ occurring outside 
$\{e\}^{r/3}$.''  $C_{ij}^A$ and $C_{ij}^B$ are not necessarily disjoint, but 
\begin{equation} \label{E:CijAB}
  C_{ij} = C_{ij}^A \cup C_{ij}^B.
\end{equation}

If $A$ (or $B$) occurs on some $\mathcal{E}
\subset \mR$ in the bottom layer of the $A$ (or $B$) pair, then by
(\ref{E:FKGcoup}), $A$ (or $B$) occurs on the same $\mathcal{E}$ in the top layer
of the same pair.  Hence
\[
  C_{00} \quad  \subset \quad C_{10} \cap C_{01} \quad  \subset \quad 
    C_{10} \cup C_{01} \quad  \subset \quad C_{11}.
\]
It follows
that $\{0,1\}^{\mathcal{U}} \times (\{0,1\}^{\mathcal{S}})^{4}$ is the disjoint
union of the sets
\[
  C_{11}^{c}, \quad C_{00}, \quad C_{10} \backslash C_{01}, \quad C_{01}
    \backslash C_{10}, \quad (C_{10} \cap C_{01}) \backslash C_{00}, \quad
    C_{11} \backslash (C_{10} \cup C_{01}).
\]
We next consider conditioning on each of these.  

If $(\zeU,\zeS^{1},\zeS^{0},\zeSt^{1},\zeSt^{0}) \in C_{11}^{c}$, then split
occurrence cannot occur no matter what layers are chosen, and
\begin{align} \label{E:C11c}
  \mathbb{P}_{\mathcal{S},e} &\bigl( C(A,B,r,\mathcal{S} \cup \{e\}) \mid
    \omU = \zeU,\omS^{1} = \zeS^{1},\omS^{0} = \zeS^{0},\omSt^{1} = \zeSt^{1},
    \omSt^{0} = \zeSt^{0} \bigr) \\
  &= 0 \notag \\
  &= \mathbb{P}_{\mathcal{S},e}\bigl( C(A,B,r,\mathcal{S}) \mid
    \omU = \zeU,\omS^{1} = \zeS^{1},\omS^{0} = \zeS^{0},\omSt^{1} = \zeSt^{1},
    \omSt^{0} = \zeSt^{0} \bigr). \notag
\end{align}
    
If $(\zeU,\zeS^{1},\zeS^{0},\zeSt^{1},\zeSt^{0}) \in C_{00}$, then split occurence
will occur regardless of what layers are chosen, and
\begin{align} \label{E:C00}
  \mathbb{P}_{\mathcal{S},e} &\bigl( C(A,B,r,\mathcal{S} \cup \{e\}) \mid
    \omU = \zeU,\omS^{1} = \zeS^{1},\omS^{0} = \zeS^{0},\omSt^{1} = \zeSt^{1},
    \omSt^{0} = \zeSt^{0} \bigr) \\
  &= 1 \notag \\
  &= \mathbb{P}_{\mathcal{S},e}\bigl( C(A,B,r,\mathcal{S}) \mid
    \omU = \zeU,\omS^{1} = \zeS^{1},\omS^{0} = \zeS^{0},\omSt^{1} = \zeSt^{1},
    \omSt^{0} = \zeSt^{0} \bigr). \notag
\end{align}

If $(\zeU,\zeS^{1},\zeS^{0},\zeSt^{1},\zeSt^{0}) \in C_{10} \backslash 
C_{01}$, then for split occurrence we must choose the top layer in the $A$ pair;
we have
\begin{align} \label{E:C10only}
  \mathbb{P}_{\mathcal{S},e} &\bigl( C(A,B,r,\mathcal{S} \cup \{e\}) \mid
    \omU = \zeU,\omS^{1} = \zeS^{1},\omS^{0} = \zeS^{0},\omSt^{1} = \zeSt^{1},
    \omSt^{0} = \zeSt^{0} \bigr) \\
  &= P(\ome = 1 \mid \omU = \zeU) \notag \\
  &= \mathbb{P}_{\mathcal{S},e}\bigl( C(A,B,r,\mathcal{S}) \mid
    \omU = \zeU,\omS^{1} = \zeS^{1},\omS^{0} = \zeS^{0},\omSt^{1} = \zeSt^{1},
    \omSt^{0} = \zeSt^{0} \bigr). \notag
\end{align}
and similarly if $(\zeU,\zeS^{1},\zeS^{0},\zeSt^{1},\zeSt^{0}) 
\in C_{01} \backslash C_{10}$, where we must choose the top layer in the $B$ pair.

If $(\zeU,\zeS^{1},\zeS^{0},\zeSt^{1},\zeSt^{0}) \in (C_{10} \cap 
C_{01}) \backslash C_{00}$, then we must choose the top layer in at least one pair,
and
\begin{align} \label{E:either}
  \mathbb{P}_{\mathcal{S},e} &\bigl( C(A,B,r,\mathcal{S} \cup \{e\}) \mid
    \omU = \zeU,\omS^{1} = \zeS^{1},\omS^{0} = \zeS^{0},\omSt^{1} = \zeSt^{1},
    \omSt^{0} = \zeSt^{0} \bigr) \\
  &= \mathbb{P}_{\mathcal{S},e}(\ome = 1 \text{ or } 
    \omet = 1 \mid \omU = \zeU) \notag \\
  &= 1 - P(\ome = 0 \mid \omU = \zeU)^{2} \notag \\
  &\geq P(\ome = 1 \mid \omU = \zeU) \notag \\
  &= \mathbb{P}_{\mathcal{S},e}\bigl( C(A,B,r,\mathcal{S}) \mid
    \omU = \zeU,\omS^{1} = \zeS^{1},\omS^{0} = \zeS^{0},\omSt^{1} = \zeSt^{1},
    \omSt^{0} = \zeSt^{0} \bigr). \notag
\end{align}

If $(\zeU,\zeS^{1},\zeS^{0},\zeSt^{1},\zeSt^{0}) \in C_{11} \backslash 
(C_{10} \cup C_{01})$, then we must choose the top layer in both pairs, so
the analog of (\ref{E:either}) fails because we would
have to replace ``or'' with ``and'' in the second line; the third line would be 
$P(\ome = 1 \mid \omU = \zeU)^{2}$ and
the inequality would then go the wrong way.  However, by (\ref{E:CijAB}), 
\begin{equation} \label{E:both}
  (\zeU,\zeS^{1},\zeS^{0},\zeSt^{1},\zeSt^{0}) \in \bigl( C_{11}^A \bs (C_{01} \cup
    C_{01}) \bigr) \cup \bigl( C_{11}^B \bs (C_{01} \cup
    C_{01}) \bigr).
\end{equation}
Combining (\ref{E:C11c})--(\ref{E:both}) we
see that
\begin{align} \label{E:compare}
  P_{\mathcal{S}}&(A \circ_{r,\mathcal{S}} B) \\
  &= \mathbb{P}_{\mathcal{S},e} \bigl( C(A,B,r,\mathcal{S}) \bigr) \notag \\
  &\leq \mathbb{P}_{\mathcal{S},e} \bigl( C(A,B,r,\mathcal{S} \cup \{e\} ) \bigr)
    + \mathbb{P}_{\mathcal{S},e} \bigl( (\omU,\omS^{1},\omS^{0},\omSt^{1},
    \omSt^{0}) \in C_{11}^A \bs (C_{10} \cup C_{01}), \ome = \omet = 1 \bigr)
    \notag \\
  &\qquad + \mathbb{P}_{\mathcal{S},e} \bigl( (\omU,\omS^{1},\omS^{0},\omSt^{1},
    \omSt^{0}) \in C_{11}^B \bs (C_{10} \cup C_{01}), \ome = \omet = 1 \bigr)
    \notag \\
  &= P_{\mathcal{S} \cup \{e\} }(A \circ_{r,\mathcal{S} \cup \{e\} } B)
    + \mathbb{P}_{\mathcal{S},e} \bigl( (\omU,\omS^{1},\omS^{0},\omSt^{1},
    \omSt^{0}) \in C_{11}^A \bs (C_{10} \cup C_{01}), \ome = \omet = 1 \bigr)
    \notag \\
  &\qquad + \mathbb{P}_{\mathcal{S},e} \bigl( (\omU,\omS^{1},\omS^{0},\omSt^{1},
    \omSt^{0}) \in C_{11}^B \bs (C_{10} \cup C_{01}), \ome = \omet = 1 \bigr).
    \notag 
\end{align}
To obtain (\ref{E:inducstep}) it is now sufficient to show that the couplings 
$\hat{P}_{\zeU}$ can be chosen so that
\begin{equation} \label{E:suffices}
  \mathbb{P}_{\mathcal{S},e} \bigl( (\omU,\omS^{1},\omS^{0},\omSt^{1},
    \omSt^{0}) \in C_{11}^A \bs (C_{10} \cup C_{01}), \ome = \omet = 1 \bigr) \leq 
    c_{23} e^{-\ep_{14} r} P_{\mathcal{S}}(A \circ_{r,\mathcal{S}} B),
\end{equation}
and similarly for $C_{11}^B$ in place of $C_{11}^A$.
By virtue of ({\ref{E:Cprobs}), (\ref{E:suffices}) says roughly that given that $A
\circ_{r,\mathcal{S}} B$ occurs in the configuration using the top layer of each
pair, with $A$ occurring far from $e$, it is exponentially unlikely that $A
\circ_{r,\mathcal{S}} B$ fails to occur (on the same separated sets of bonds $\mE$
and $\mF$ for $A$ and $B$ respectively, actually) when the top layer is replaced by
the bottom layer in the $A$ pair.  This, we will see, is because the top and
bottom layers are likely equal far from $e$. 

For the proof of (\ref{E:suffices}) one key is the next proposition.  The idea
is as follows.  We wish to consider the effect of the
configuration in a region
$\mF$ on probabilities of events occurring on a distant region $\mE$, with
$\mE,\mF$ contained in some larger region $\mR$.  In ordinary (unsplit)
configurations, this effect is exponentially small provided the ratio strong mixing
property holds.  Suppose, though, that we have split the configuration on a subset
$\mS$ of $\mR$, and suppose that $\mE$ consists of unsplit bonds.  We then have
\emph{two} configurations on the split portion of $\mF$ exerting their influence
on probabilities for events occurring on $\mE$, and it is not \emph{a priori}
clear under ratio strong mixing that this influence is still exponentially small. 
The proposition guarantees this smallness, at least when the influence is
measured additively, not using the ``ratio'' form of influence. 

\begin{proposition} \label{P:splitmix}
  Let $\mkR, \mkS, \mkM$ be as in the standard bond percolation setup and
  let $\mkN$ be a neighborhood collection for $\mkS$.
  Suppose each measure in $\mkM$ satisfies the FKG lattice condition,
  each measure in $\mkM \cup \mkM^+(\mkS,\mkR)$
  has the Markov property for blocking sets, and 
  $\mkM_w^+(\mkN,\mkR)$ has uniform exponential decay of connectivity.  There exist
  $c_{24}, \ep_{15}$ as follows.  Let $\mR \in \mkR, P \in \mkM_{\mR}$ and $\mS \in
  \mkS_{\mR}$, and let $\mT = \mR \bs \mS$.
  Suppose that or some
  $c>1$, for all $r>0$ and $e \in \mT$, $e$ is neighborhood-appendable to $\mS$ for
  $(c,r,\mR,P)$.
  Let $\mE \subset \mT, \mF \subset \mR$  and $G \in \mG_{\mE}$.  
  Then for every choice of configurations 
  $\rho_{\mF \cap \mT}, \rho_{\mF \cap \mS}, \tilde{\rho}_{\mF \cap \mS}$,
  \begin{align}
    |P_{\mS}&(G \mid \omega_{\mF \cap \mT} = \rho_{\mF \cap \mT},
      \omega_{\mF \cap \mS} = \rho_{\mF \cap \mS}, 
      \tilde{\omega}_{\mF \cap \mS} = \tilde{\rho}_{\mF \cap \mS}) - 
      P(G)| \notag \\
    &\leq c_{24} \sum_{x \in V(\mE), y \in V(\mF)} e^{-\ep_{15} d_{\mR}(x,y)}.
      \notag 
  \end{align}
\end{proposition}

For the proof we need the following.

\begin{proposition} \label{P:strongmix}
  Let $\mkR, \mkS, \mkM$ be as in the standard bond percolation setup and
  let $\mkN$ be a neighborhood collection for $\mkS$.
  Suppose each measure in $\mkM$ satisfies the FKG lattice condition,
  each measure in $\mkM^+(\mkS,\mkR)$
  has the Markov property for blocking sets, and 
  $\mkM_w^+(\mkN,\mkR)$ has uniform exponential decay of connectivity.  Then 
  $\mkM^+(\mkS,\mkR)$
  has the ratio strong mixing property.
\end{proposition}
\begin{proof}
This is proved in (\cite{Al01fin}, Theorem 1.6) in the special case of the FK
model with site boundary conditions, with $\mkS_{\mR} = \{ \mR \}$.
In that special case, not all measures in
$\mkM$ satisfy the FKG lattice condition, and the following property of the FK
model under site boundary conditions is implicitly used instead:  for every $\mR
\in \mkR$, every $P
\in \mkM_{\mR}$, every $\mS \subset \mR$ and every configuration $\rho_{\mR \bs
\mS}$ with $P(\omega_{\mR \bs \mS} = \rho_{\mR \bs
\mS}) > 0$,  the measure $P(\cdot \mid \omega_{\mR \bs \mS} = \rho_{\mR \bs
\mS})$ is FKG-dominated by the wired-boundary measure conditioned on 
$\omega_{\mR \bs \mS} = \rho_{\mR \bs \mS}$, which does satisfy the FKG lattice
condition.  The arguments used in
\cite{Al01fin}, including those from
\cite{Al98mix} cited in \cite{Al01fin}, are essentailly unchanged under the
assumptions of the present proposition.
\end{proof}

\begin{remark} \label{R:QnotinM}
The proof of Proposition \ref{P:strongmix}, as given in \cite{Al01fin} and
\cite{Al98mix}, shows that under the hypotheses given, the
ratio weak mixing statement (\ref{E:rstrongmix}) actually holds for measures $Q
\notin \mkM^+(\mkS,\mkR)$, of form $Q = P(\omega_{\mS \cup \mE} 
\in \cdot \mid \omega_{\mR
\bs (\mS \cup \mE)} = \rho_{\mR \bs (\mS \cup \mE)}^i)$ for some
$\mR \in \mkR, \mS \in \mkS_{\mR}, \mE \subset \mR$ and $P \in \mkM_{\mR}$,
so long as $Q$ has the
Markov property for sets blocking $\mE$.
That is, if we consider 
configurations on the region $\mS$ and view $\mE \bs \mS$ as part
of the ``partial boundary'' $\mR \bs \mS$ of $\mS$, then we
can limit the influence of both the boundary and non-boundary portions ($\mE \bs
\mS$ and $\mE \cap \mS)$ of $\mE$ on distant events, so long as the influence of
$\mE$ can be blocked by a barrier of closed bonds.  Further, in this situation, we
need not assume that all of $\mkM^+(\mkS,\mkR)$ has the Markov property for
blocking sets; the Markov property for $Q$ for sets blocking $\mE$ is sufficient.
For example, if $\mkS_{\mR}$
consists of all SLC subsets of $\mR$ then it is not necessary that $\mS \cup \mE$
be SLC; instead it suffices that $\mS$ be SLC, so long as $Q$ has the Markov
property for sets blocking
$\mE$.  See Figure \ref{F:markov}.
\end{remark}

\begin{proof}[Proof of Proposition \ref{P:splitmix}.]
First observe that $P_{\mS}(G) = P(G)$, since $G$ occurs on unsplit bonds only.  
Also, for fixed $e$ the measure $P_{\mS}(\ome \in
\cdot)$ is FKG-dominated by $P_{\mS}(\ome \in \cdot \mid \omega_{\mF \cap \mT} =
\rho_{\mF
\cap \mT}^1, \omega_{\mF \cap \mS} = \rho_{\mF \cap \mS}^1, 
\tilde{\omega}_{\mF \cap \mS} = \tilde{\rho}_{\mF \cap \mS}^1)$, and it
FKG-dominates $P_{\mS}(\ome \in \cdot \mid \omega_{\mF \cap \mT} =
\rho_{\mF
\cap \mT}^0, \omega_{\mF \cap \mS} = \rho_{\mF \cap \mS}^0, 
\tilde{\omega}_{\mF \cap \mS} = \tilde{\rho}_{\mF \cap \mS}^0)$.  It follows that
there exists a 3-way FKG coupling of these measures, that is, a coupling in
which the configuration under $P_{\mS}(\ome \in
\cdot)$ is always sandwiched between the other two configurations, in the usual
partial ordering of configurations.  A similar 3-way ``sandwiching'' coupling can
be created using $P_{\mS}(\ome \in \cdot \mid \omega_{\mF \cap \mT} = \rho_{\mF
\cap \mT},
\omega_{\mF \cap \mS} = \rho_{\mF \cap \mS}, 
\tilde{\omega}_{\mF \cap \mS} = \tilde{\rho}_{\mF \cap \mS})$ in place of 
$P_{\mS}(\ome \in \cdot)$.  
It follows easily from the existence of these sandwiching couplings that 
\begin{align} \label{E:onebond}
    \bigl|P_{\mS}&(G \mid \omega_{\mF \cap \mT} = \rho_{\mF \cap \mT},
      \omega_{\mF \cap \mS} = \rho_{\mF \cap \mS}, 
      \tilde{\omega}_{\mF \cap \mS} = \tilde{\rho}_{\mF \cap \mS}) - 
      P_{\mS}(G) \bigr| \\
    &\leq \sum_{e \in \mE} \Bigl(
      P_{\mS}(\omega_e = 1 \mid \omega_{\mF \cap \mT} 
      = \rho_{\mF \cap \mT}^1,
      \omega_{\mF \cap \mS} = \rho_{\mF \cap \mS}^1, 
      \tilde{\omega}_{\mF \cap \mS} = \tilde{\rho}_{\mF \cap \mS}^1) \notag \\
    &\qquad \qquad - P_{\mS}(\omega_e = 1 \mid \omega_{\mF \cap \mT} 
      = \rho_{\mF \cap \mT}^0,
      \omega_{\mF \cap \mS} = \rho_{\mF \cap \mS}^0, 
      \tilde{\omega}_{\mF \cap \mS} = \tilde{\rho}_{\mF \cap \mS}^0) 
      \Bigr). \notag 
\end{align}
Thus we may assume $\mE$ consists of a single bond $e \in \mT$ and $G = [\ome =
1]$.  Define
$\mU = \mT \bs \{e\}$ and $r = d_{\mR}(e,\mF)$.
Since by assumption $e$ is neighborhood-appendable to $\mS$, there exist $c_{25}>1$
and a $c_{25}$-approximate
$r/4$-neighborhood $\mW$ of $e$ in $\mR$ such that either $\mW$ or $\mW \bs \mS$ is
blockable in $\mS \cup \mW$ under $P$.  By enlarging $\mF$ if
necessarily, we may assume that
\begin{equation} \label{E:newF}
  \mF = (\mU \bs \mW) \cup (\mS \bs \mB_{\mR}(e,r));
\end{equation}
see Figure \ref{F:Fchoice}.  Then
\begin{equation} \label{E:separation}
  d_{\mR}(\mT \bs \mF, \mF \cap \mS) > \frac{r}{2}.
\end{equation}
Write $P^i_{\mS}$ for $P_{\mS}(\cdot \mid \omega_{\mF
\cap \mT}  = \rho_{\mF \cap \mT}^i,
\omega_{\mF \cap \mS} = \rho_{\mF \cap \mS}^i, 
\tilde{\omega}_{\mF \cap \mS} = \tilde{\rho}_{\mF \cap \mS}^i), i = 0,1$.  
To obtain the desired bound on (\ref{E:onebond}),
it is thus sufficient to show that 
\begin{equation} \label{E:suffcond}
  |P_{\mS}^1(\ome = 1) - P_{\mS}^0(\ome = 1)| \leq 
    c_{26} \sum_{x \in V(e),y \in V(\mF)} e^{-\ep_{15} d_{\mR}(x,y)}.
\end{equation}

\input{epsf}
\begin{figure} 
\epsfxsize=3.5in
\begin{center}
\leavevmode
\epsffile{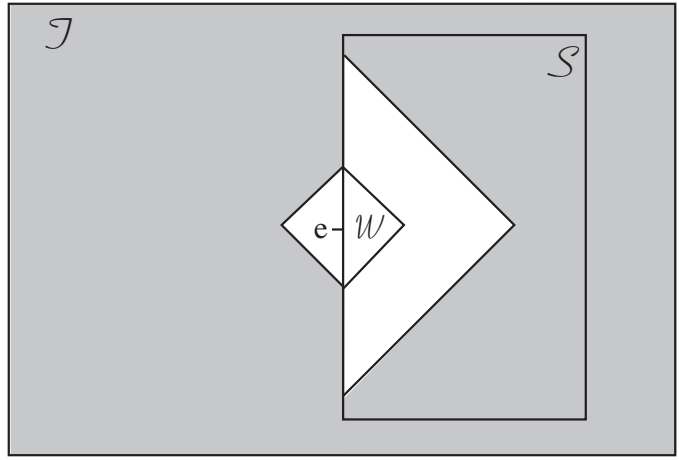}
\end{center}
\caption{$\mF$ is the shaded region.  $\mS$ is the inner rectangle, and $\mT$ is
the complement of $\mS$ in the full rectangle $\mR$.  The large triangular region
is $\mS \cap B_{\mR}(e,r)$.  The figure is not precisely to scale.}
\label{F:Fchoice}
\end{figure}

Define
\[
  g(\rho_e,\rho_{\mU \bs \mF}) = P(\omega_{\mF \cap \mS} = \rho_{\mF
    \cap \mS}^1
    \mid \omega_{\mF \cap \mU} = \rho_{\mF \cap \mU}^1, \omega_{\mU \bs \mF}
    = \rho_{\mU \bs \mF},
    \ome = \rho_e).
\]
Then by the nature of the $\mS$-split measure,
\begin{equation} \label{E:useg}
  P_{\mS}^1(\ome = 1) = \frac{E(\delta_{[\ome = 1]}g 
    \mid \omega_{\mF} =
    \rho_{\mF}^1)}{E(g \mid \omega_{\mF} = \rho_{\mF}^1)},
\end{equation}
where $E$ denotes expectation with respect to $P$.  Let 
$P_1 = P(\omega_{\mR \bs (\mF \cap \mU)} \in \cdot \mid \omega_{\mF
\cap \mU} =\rho_{\mF \cap \mU}^1)$, so that $g(\rho_e,\rho_{\mU \bs \mF}) = 
P_1(\omega_{\mF \cap \mS} = \rho_{\mF \cap \mS}^1 \mid \omega_{\mU \bs \mF}
= \rho_{\mU \bs \mF}, \ome = \rho_e)$. 
By Proposition \ref{P:strongmix}, $\mkM^+(\mkS,\mkR)$ has the
ratio strong mixing property, and we would like to apply this property to control
the effect of this conditioning of $P_1$, but the location $\mR \bs (\mF \cap
\mU) = \mS \cup \mW$ of the configurations under
$P_1$ need not be in $\mkS_{\mR}$.  However, by Remark
\ref{R:QnotinM} (with $\mW$ or $\mW \bs \mS$, whichever is blockable under the
admissibility assumption, in place of $\mE$) this is not a problem.  Thus we have
for all $\rho_e,\rho_{\mT \bs \mF}$,
\[
  \left| \frac{g(\rho_e,\rho_{\mT \bs \mF})}{P_1(\omega_{\mF \cap \mS} = 
    \rho_{\mF \cap \mS}^1)} - 1 \right| \leq
    c_{27} \sum_{x \in V(\mT \bs \mF),y \in V(\mF \cap \mS)} e^{-\lambda
    d_{\mR}(x,y)} \leq c_{28}e^{-\ep_{16}r},
\]
so that $g$ is nearly a constant in (\ref{E:useg}), and we obtain
\[
  \left| \frac{P_{\mS}^1(\ome = 1)}{P(\ome = 1 \mid \omega_{\mF} = \rho_{\mF}^1)}
    - 1 \right| \leq c_{29}e^{-\ep_{17}r}.
\]
Similarly,
\[
  \left| \frac{P_{\mS}^0(\ome = 1)}{P(\ome = 1 \mid \omega_{\mF} = \rho_{\mF}^0)}
    - 1 \right| \leq c_{29}e^{-\ep_{17}r}.
\]
Next we apply Proposition \ref{P:strongmix} to the collection $\mkS_0 = \{ \{ \mR
\}: \mR \in \mkR \}$, noting that $\mkM^+(\mkS_0,\mkR) = \mkM$, to obtain
\[
  \left| \frac{P(\ome = 1 \mid \omega_{\mF} = \rho_{\mF}^0)}{P(\ome = 1 \mid
    \omega_{\mF} = \rho_{\mF}^1)}
    - 1 \right| \leq c_{30}e^{-\ep_{18}r}.
\]
Combining the last 3 inequalities we obtain
\[
  |P_{\mS}^1(\ome = 1) - P_{\mS}^0(\ome = 1)| \leq c_{31}e^{-\ep_{19}r},
\]
which establishes (\ref{E:newF}).
\end{proof}

We now turn to the proof of (\ref{E:suffices}).  Since $\mkM_w^+(\mkN,\mkR)$ has
uniform exponential decay of connectivity, it follows easily using the FKG
property that if we take the class of balls $\mkB = \mkB(\mkR) = \{
B_{\mR}(x,s): \mR \in \mkR, x
\in V(\mR), s > 1 \}$ in place of $\mkN$, the class $\mkM^+(\mkB,\mkR)$ has uniform
exponential decay as well, that is, there exists
$C, \lambda > 0$ such that 
\[
  P(x \lra y \text{ via a path in } \mB \mid \omega_{\mR \bs \mB} = \rho_{\mR \bs
    \mB})  \leq C e^{-\lambda d_{\mR}(x,y)}
\]
for all $\mR \in \mkR$, all $\mB \subset \mR$ in $\mkB$, all $x,y \in V(\mB)$ and
all boundary conditions
$\rho$.  

Recall that for (\ref{E:suffices}) we have a fixed $\mR \in \mkR$, a filling
sequence $\mV = \{ b_1,..,b_n \}$ and $\mS = \{ b_1,..,b_{k-1} \}, e = b_k$. 
Let $m = \lfloor r/24 \rfloor$.  For
$x \in V(\mR)$ let $\mQ_x = B_{\mR}(x,3m) \cap \mT$ and $\mA_x = \mR \bs
B_{\mR}(x,3m)$, and
for each configuration $\rho_{\mQ_x}$ let
\[
  \phi_x(\rho_{\mQ_x}) 
  = P\left( x \lra y \text{ for some } y \in V(\mR) \text{ with }
    d_{\mR}(x,y) \geq m \ \Big| \ \omega_{\mQ_x} = \rho_{\mQ_x}, 
    \omega_{\mA_x} = \rho_{\mA_x}^1 \right). 
\]
Then since $B_{\mR}(x,3m) \in \mkB$,
\begin{align} \label{E:phixbound}
  E_{\mR,w}(\phi_x(\omega_{\mQ_x})) &\leq
    P\left( x \lra y \text{ for some } y \in V(\mR) \text{ with }
    d_{\mR}(x,y) \geq m \ \Big| \  
    \omega_{\mA_x} = \rho_{\mA_x}^1 \right) \\
  &\leq c_{32}m^{d-1}\ Ce^{-\lambda m}.  \notag 
\end{align}
We say that $\mQ_x$ is \emph{connection--inducing} in a configuration
$\rho_{\mQ_x}$ if $\phi_x(\rho_{\mQ_x}) \geq C e^{-\lambda m/2}$.
Roughly speaking, $\mQ_x$ is connection-inducing if $\rho_{\mQ_x}$ either contains
a long open path starting in $Q_{m}(x)$, or contains enough segments of such a
path that the conditional probability for such a path given $\rho_{\mQ_x}$ is
greatly increased above the unconditional probability.  If $\mQ_x$ is not
connection--inducing we say it is \emph{insulating}.  Note that $\phi_x$
is an increasing function.  Hence using the FKG inequality and (\ref{E:phixbound}),
for every $x$ and $\zeta$,
\begin{align} \label{E:Markovineq}
  P_{\mR,\zeta}(\mQ_x \text{ is connection--inducing}) 
    &= P_{\mR,\zeta}(\phi_x(\omega_{\mQ_x}) \geq C e^{-\lambda m/2}) \\ 
  &\leq \frac{E_{\mR,w}(\phi_x(\omega_{\mQ_x}))}{C e^{-\lambda m/2}} \notag \\
  &\leq c_{33}e^{-\lambda m/3}. \notag
\end{align}

The idea can now be sketched as follows.  Recall that $\mU = \mT \bs \{e\}$. 
Suppose
$A$ and $B$ occur $(\mS \cup \{e\})$-split at separation
$r$ in $(\omega_{\mU},(1,\omega_{\mS}^1),(1,\tilde{\omega}_{\mS}^1))$ with $A$
occurring far from $e$, that is, 
$(\omega_{\mU},\omega_{\mS}^1,\omega_{\mS}^0,\tilde{\omega}_{\mS}^1,
\tilde{\omega}_{\mS}^0) \in C_{11}^A$, with $A$ occuring on some set $\mE$ and $B$
on some set $\mF$, with $d_{\mR}(e,\mE) \geq r/3$.  Since $P$ has the FKG property
and the Markov property for blocking sets, the
FKG couplings $\hat{P}_{\zeta_{\mU}}, \zeta_{\mU} \in \{0,1\}^{\mU}$,
of the two layers of the $A$ pair can be
chosen so that $\omega_{\mS}^1 = \omega_{\mS}^0$ outside
$C_e = C_e((\zeta_{\mU},1,\omega_{\mS}^1))$, which is the cluster of $e$ for the
top layer together with the unsplit bonds.  (The construction of such couplings is
standard--see e.g. \cite{AC}, \cite{Ne}.)  Suppose
now that the $A$ pair fails to couple on $\mE$, that is, $\omega_{\mE \cap \mS}^1
\neq \omega_{\mE \cap \mS}^0$.  This means that the cluster $C_e$ intersects
$\mE$, i.e. there is an open path from $e$ to $\mE$ in the top layer of
the $A$ pair.  Since $d(\mE,\mF)
\geq r$, a segment of this path is far from $\mE, \mF$ and $e$.  This segment
may be partly in $\mU$ and partly in $\mS$.  There are two possibilities: either
the segment is substantially helped to exist by some connection-inducing region
of $\mU$, or the portion in $\mS$ exists without such help, that is, the
relevant regions of $\mU$ are insulating.   The first is
exponentially unlikely, even conditionally on the occurrence of $A \circ_{r,\mS}
B$ on $\mE \cup \mF$ in the top layers of the $A$ and $B$ pairs, by Proposition
\ref{P:splitmix}.  The second is also (conditionally) exponentially unlikely, by
the definition of insulating and the ratio strong mixing property.

Turning to the details, we let $D_{11}^A$ denote the set of all $(\omega_{\mU},
\omega_{\mS}^1,\omega_{\mS}^0,
\tilde{\omega}_{\mS}^1, \tilde{\omega}_{\mS}^0)$ such that there exist $\mE,\mF$
for which $A$ and $B$ occur $(\mS \cup \{ e \})$-split at separation $r$
(which is at least $24m$) in
$(\omega_{\mU},(1,\omega_{\mS}^1),(1,\tilde{\omega}_{\mS}^1))$, i.e. using the
two top layers, with
$A$ occuring on $\mE$, $B$ occurring on $\mF$, $d_{\mR}(e,\mE) > 8m$, and
$\omega_{\mE}^1 \neq
\omega_{\mE}^0$, and for some $z \in V(\mD^r(\mR))$, we have that $\mQ_z$ is
connection-inducing, $d_{\mR}(z,e) > 3m$ and $d_{\mR}(z,\mE \cup \mF) > 4m$.   Let
$E_{11}^A$ be the set of all
$(\omega_{\mU},
\omega_{\mS}^1,\omega_{\mS}^0,
\tilde{\omega}_{\mS}^1, \tilde{\omega}_{\mS}^0)$ such that there exist $\mE,\mF$
for which $A$ and $B$ occur $(\mS \cup \{ e \})$-split at separation $r$ in
$(\omega_{\mU},(1,\omega_{\mS}^1),(1,\tilde{\omega}_{\mS}^1))$ with
$A$ occuring on $\mE$, $B$ occurring on $\mF$, $d_{\mR}(e,\mE) > 8m$, and
$\omega_{\mE}^1
\neq \omega_{\mE}^0$, but for no choice of such $\mE,\mF$ does there exist $z$ as
above.  $D_{11}^B$ and $E_{11}^B$ are defined analogously.
Then 
\[
  C_{11}^A \bs (C_{01} \cup C_{10}) \subset D_{11}^A \cup E_{11}^A,
\]
and similarly with $B$in place of $A$.  
Therefore provided $c_{7}$ is sufficiently large, by Proposition
\ref{P:splitmix} and (\ref{E:Markovineq}),
\begin{align} \label{E:D11bound}
  \mathbb{P}_{\mS,e}&((\omU,\omS^{1},\omS^{0},\omSt^{1},\omSt^{0}) \in D_{11}^A,
    \ome = \omet = 1) \\
  &\leq \sum_{z \in V(\mD^r(\mR)): d_{\mR}(z,e) > 3m}
    \mathbb{P}_{\mS,e}\Bigl(\ome = \omet = 1, \mQ_z \text{ is
    connection-inducing in } (\omega_{\mU},1), \notag \\
  &\qquad \qquad \qquad \qquad \qquad \qquad A \circ_{r,\mS} B \text{ occurs on } 
    B_{\mR}(z,4m)^c) \text{ in }
    (\omega_{\mU},(1,\omega_{\mS}^1),(1,\tilde{\omega}_{\mS}^1)) \Bigr) \notag \\
  &\leq \sum_{z \in V(\mV)} P_{\mS}(A \circ_{r,\mS} B) \cdot \Bigl( 
    \mathbb{P}_{\mS,e} \bigl( \mQ_z \text{ is connection-inducing }
    \bigr) \notag \\
  &\qquad \qquad \qquad \qquad \qquad \qquad \qquad \qquad
    + C \sum_{x \in V(B_{\mR}(z,3m)), y \in V(B_{\mR}(z,4m)^c))} e^{-\lambda
    d_{\mR}(x,y)} \Bigr) \notag \\
  &\leq c_{34}|\mV| e^{-\ep_{20}r} P_{\mS}(A \circ_{r,\mS} B) \notag \\
  &\leq c_{35} e^{-\ep_{21}r} P_{\mS}(A \circ_{r,\mS} B). \notag
\end{align}
Here we have used the fact that when $d_{\mR}(z,e) > 3m$ we have $e \notin \mQ_z$,
so whether or not
$\mQ_z$ in connection-inducing in $(\omU,\ome)$ does not depend on $\ome$.  Next,
let $H_z$ denote the event that $z \lra y
\text{ in } (\omega_{\mU},1, \omega_{\mS}^1)$
for some $y$ with $|y - z| \geq m$, $I_z$ the event that $\mQ_z$ is insulating in
$(\omega_{\mU},1)$ and $J_z$ the event that $A \circ_{r,\mS} 
B$ occurs on $B_{\mR}(z,4m)^c$ in
$(\omega_{\mU},(1,\omega_{\mS}^1),(1,\tilde{\omega}_{\mS}^1))$.  
Observe that if $(\omU,\omS^{1},\omS^{0},\omSt^{1},\omSt^{0}) \in E_{11}^A$ then
there is an open path from $e$ to $\mE$ in the top layer of the $A$ pair, and this
path must path through a site $z$ at $d_{\mR}$-distance approximately
$4m$ from $\mE$; this forces the event $H_z \cap I_z \cap J_z$ to occur. 
(Here we use the fact that $\omega_{\mS}^1 = \omega_{\mS}^0$ outside $C_e$.)  
Hence using the conditional independence inherent in the $\mS$-split structure,
\begin{align} \label{E:E11bound}
  \mathbb{P}_{\mS,e}&((\omU,\omS^{1},\omS^{0},\omSt^{1},\omSt^{0}) \in E_{11}^A,
    \ome = \omet = 1) \\
  &\leq \sum_{z \in V(\mD^r(\mR)): d_{\mR}(z,e) > 3m }
    \mathbb{P}_{\mS,e}( H_z \cap I_z \cap J_z, \ome = 1 ) \notag \\
  &\leq \sum_{z \in V(\mD^r(\mR)): d_{\mR}(z,e) > 3m } 
    \ \sum_{\zeta_{\mU} \in I_z} 
    \mathbb{P}_{\mS,e}(J_z, \ome = 1, \omega_{\mU} = \zeta_{\mU}) 
    \mathbb{P}_{\mS,e}(H_z \mid J_z, \ome = 1, \omega_{\mU} = \zeta_{\mU}). \notag
\end{align}
From the definition of insulating and the FKG property we have
\[
  \mathbb{P}_{\mS,e}(H_z \mid J_z, \ome = 1, \omega_{\mU} = \zeta_{\mU}) < 
    Ce^{-\lambda m/2},
\]
so (\ref{E:E11bound}) yields
\begin{align} \label{E:E11bound2}
  \mathbb{P}_{\mS,e}&((\omU,\omS^{1},\omS^{0},\omSt^{1},\omSt^{0}) \in E_{11}^A,
    \ome = \omet = 1) \\
  &\leq \sum_{z \in V(\mV) } Ce^{-\lambda m/2}\ \mathbb{P}_{\mS,e}(J_z, \ome
    = 1) \notag \\
  &\leq c_{36} |\mV| e^{-\lambda m/2} P_{\mS}(A \circ_{r,\mS} B) \notag \\
  &\leq c_{37} e^{-\ep_{22}r} P_{\mS}(A \circ_{r,\mS} B). \notag
\end{align}
Now (\ref{E:D11bound}) and (\ref{E:E11bound2}) show that 
\[
  \mathbb{P}_{\mS,e}\Bigl( (\omU,\omS^{1},\omS^{0},\omSt^{1},\omSt^{0}) \in
    C_{11}^A \bs (C_{01} \cup C_{10}), \ome = \omet = 1 \Bigr) \leq c_{38}
    e^{-\ep_{23}r} P_{\mS}(A \circ_{r,\mS} B)
\]
and a similar bound holds for $C_{11}^B$, so we have established
(\ref{E:suffices}), and thus also (\ref{E:inducstep}) and then (\ref{E:induc}). 
This proves Theorem
\ref{T:finitevol}.

\begin{proof}[Proof of Theorem \ref{T:FKZd}.]
Let $\mkR,\mkS,\mkV,\mkM,\mkN$ be as in Example
\ref{Ex:slcprops}.  

Suppose first that $\rho$ is wired or all external fields are 0.
As noted in 
Example \ref{Ex:slcprops}, $(\mkR,\mkS,\mkV,$ $\mkM,\mkN)$ is then
filling-compatible.  All finite-volume measures $P_{\mB,\rho}$,
with $\mB$ finite and  $\rho$ a bond boundary condition, satisfy the FKG
lattice condition (\cite{Fo2}; see \cite{ACCN}).  It follows that $P$ has uniform
exponential decay of connectivity for the class $\overline{\mkR}$ with arbitrary
bond boundary conditions, not just wired.  By Lemma
\ref{L:markov}, since the set $\mR \bs \mS$ is connected and abuts $\mR^c$ for all
$\mR \in \mkR$ and $\mS \in \mkS_{\mR}$, every measure in $\mkM\ \cup\
\mkM^+(\mkS,\mkR)$ has the Markov property for blocking sets.  Thus in this case
(i) follows from Theorem
\ref{T:finitevol}.

For the remaining case, suppose $\rho$ is free and not all external fields are 0. 
Consider a measure $Q = P_{\mR,f}(\omega_{\mS} \in \cdot \mid \omega_{\mR \bs \mS}
= \rho_{\mR \bs \mS}^1) \in \mkM^+(\mkS,\mkR)$.  As we have noted, $\mR \bs \mS$
is connected so the effective boundary condition on $\mS$ given by $Q$ is a
unique-cluster one. Let $(\mX, \mY, \mZ)$ be a blocking partition of $\mS$, and
suppose $\omega = 0$ on $\mY$.  We use the notation of the proof of Lemma
\ref{L:markov}.  Note that the set of bonds in $C_u$ is $\mR \bs \mS$.  The 
factoring of the weight of a cluster
$\hat{C}$ described in that proof is not necessarily valid; the clusters
$C_i$ effectively interact via the value $s(\hat{C})$.  More precisely,
conditionally on
$\{ C_j, j \neq i \}$, the effective weight attached to $C_i$ depends on $s(C_u) +
\sum_{j \neq i} s(C_j)$; since some $C_m$'s may be in $\mX$ and others in $\mZ$,
this means the Markov property for blocking sets need not hold.   However, letting
$h_k$ be the largest strictly negative external field, the effective weight of
$C_i$ is always between 1 and
$1 + O(e^{-c|h_k|s(C_u)})$, that is, the maximum influence of $\{
C_j, j \neq i \}$ on $C_i$ is exponentially small in $s(C_u)$.  Roughly speaking,
we have two situations.  If $s(C_u)$ is small relative to $r$ then since
$\diam_{\mR}(C_u) \leq s(C_u)$ the interaction between clusters $C_i$ of $\omega$
only occurs over length scales which are small relative to $r$.  If $s(C_u)$ is of
order $r$ or greater, then the above-mentioned maximum influence of $\{
C_j, j \neq i \}$ on $C_i$ is exponentially
small in $r$.  Either way, though we do not fully have the Markov property for
blocking sets, the proofs of Proposition \ref{P:splitmix}, Example
\ref{Ex:slcprops} and then Theorem \ref{T:finitevol} go through; we omit the
full details.  See (\cite{Al01fin}, Lemma 2.11(iii)) for a similar result.  Thus
(i) is proved in all cases.

For (ii) let $\mB_n = \mB([-n,n]^d)$ and let $A_n$ be the event that $A$ occurs on
$\mB_n$.  By the uniform exponential decay assumption, $P$ is the unique
infinite-volume random cluster measure at $(p,q,\{h_i\})$.  Therefore by (i),
\begin{align}
  P(A_n \circ_r B) &= \lim_m P_{\mB_m,w}(A_n \circ_r B) \notag \\
  &\leq (1 + c_{9}e^{-\ep_{5}r}) \lim_m P_{\mB_m,w}(A_n)
    P_{\mB_m,w}(B) \notag \\
  & = (1 + c_{9}e^{-\ep_{5}r})P(A_n)P(B). \notag
\end{align}
Since $A$ is locally-occurring we can now take a limit as $n \to \infty$ to obtain
(ii).
\end{proof}

\begin{proof}[Proof of Theorem \ref{T:FKZ2}.]
Let $\mkS_{\mR},\mkN$ and $\mkM$ be as in Example \ref{Ex:circbdd}.  Since $d=2$,
uniform exponential decay of connectivity for the class $\mkM^+(\mkS,\mkR)$ with
wired boundary conditions follows from the assumed infinite-volume exponential
decay \cite{Al01fin}.  As we have noted, all finite-volume measures $P_{\mB,\rho}$,
with $\mB$ finite and  $\rho$ a bond boundary condition, satisfy the FKG
lattice condition (\cite{Fo2}; see \cite{ACCN}).  We consider the case in which
the boundary is wired ($i = 1$) or there are no external fields; the case of free
boundary with external fields can be handled as in the proof of Theorem
\ref{T:FKZd}.  Example
\ref{Ex:circbdd} establishes filling-compatibility of
$(\mkR,\mkV,\mkS,\mkM,\mkN)$.  The Markov property for blocking sets for
$\mkM \cup \mkM^+(\mkS,\mkN)$ follows from Lemma \ref{L:markov}.  The theorem now
follows from Theorem \ref{T:finitevol}.
\end{proof}

We write $\Int(\gamma)$ and $\Ext(\gamma)$ for the interior and exterior of a
simple closed curve $\gamma$ in the plane.

\begin{proof}[Proof of Theorem \ref{T:farbdry}.]
For a configuration $\omega$ on $\mR$ we define the \emph{boundary cluster} 
$C_{\partial}(\omega)$ to be the set of open bonds which are connected to
$V(\mR^c)$ by a path of open bonds.  (This is a mild abuse of terminology since
$C_{\partial}$ does not necessarily consist of a single connected cluster.)  Then
$(\partial C_{\partial})^*$ includes  a finite collection of open dual circuits;
these circuits have disjoint interiors.  Let $\mkD = \mkD(\omega)$ be the set of
dual circuits in this collection which contain bonds of
$\mD^*$.  Let $\mI = \mI(\mkD(\omega)) = \cup_{\gamma \in \mkD} \
\mB(\Int(\gamma)), \mJ = \mJ(\mkD(\omega)) = \cap_{\gamma \in \mkD} \
\mB(\Ext(\gamma))$.  Conditionally on
$[\mkD = \mathbb{D}]$ for some $\mathbb{D}$, the configuration on
$\mI(\mathbb{D})$ is a free-boundary FK configuration.  Let $I_{AB}$ denote the
event that $A \circ_r B$ occurs with $B$ occuring on $\mI$ at distance $r/4$ or
more from $\mI^c$,  that is, there exist $\mE
\subset \mR, \mF \subset \mI \cap \mD$ with $d(\mE,\mF) \geq r, d(\mF,\mI^c) \geq
r/4$ such that
$A$ occurs on $\mE$ and $B$ occurs on $\mF$.  
Let $G_x$ be the event that there is an open path in
$\mR$ from $x$ to
$y$ for some $y \in V(\mR)$ with $d_{\mR}(x,y) \geq r/4$.  

Suppose $\omega \in (A \circ_r B) \bs I_{AB}$, with $A,B$ occuring on $\mE,\mF$
respectively.  Then there is an open path from the boundary $V(\mR^c) \cap V(\mR)$
to $\mF^{r/4}$, and since $d(\mE,\mF) \geq r$ a portion of length
$r/4$ of this path must be separated from both $\mE$ and $\mF$ by a distance of
more than $r/4$.  More precisely, for some $x \in V(\mD^r)$ on this path at
distance approximately $r/2$ from $\mF$ we have
$\omega
\in G_x
\cap H_x$, where $H_x$ is the event that $A \circ_r B$ occurs on $\mR \bs
B_{\mR}(x,3r/8)$.  Therefore using the FKG property,
\begin{align} \label{E:noI}
  P_{\mR,\rho} \bigl( (A \circ_r B) \bs I_{AB} \bigr) &\leq \sum_{x \in V(\mD^r)}
    P_{\mR,\rho}(G_x \mid H_x) P_{\mR,\rho}(H_x) \\
  &\leq P_{\mR,\rho}(A \circ_r B) \sum_{x \in V(\mD^r)} P_{B_{\mR}(x,3r/8),w}(G_x).
\notag
\end{align}
Since $P$ has exponential decay of connectivity, it has uniform exponential decay
for the class of all SLC subsets \cite{Al01fin} and therefore 
\[
  P_{B_{\mR}(x,3r/8),w}(G_x) \leq c_{39}e^{-\ep_{24}r}.
\]
Thus provided $c_{11}$ is large enough we have
\begin{equation} \label{E:noI2}
  P_{\mR,\rho} \bigl( (A \circ_r B) \bs I_{AB} \bigr) \leq c_{40}e^{-\ep_{25}r}
    P_{\mR,\rho}(A \circ_r B).
\end{equation}

Next we bound $P_{\mR,\rho}(I_{AB})$.
Given $\mB \subset \mR$ and a configuration
$\zeta_{\mR \bs\mB}$ on $\mR \bs \mB$, let 
\[
  A(\zeta_{\mR \bs \mB}) = \{ \zeta_{\mB}:
  (\zeta_{\mR \bs \mB}\zeta_{\mB}) \in A \}.
\]
Next let $\mI_{r/4}(\mathbb{D}) = \{ b \in \mI(\mathbb{D}), d(b,\mathbb{D}) \geq
r/4 \}$ and let $B(\mathbb{D})$ be the event that $B$ occurs on
$\mI_{r/4}(\mathbb{D})$.  We have
\begin{align} \label{E:yesI}
  P_{\mR,\rho}(I_{AB}) &\leq \sum_{\mathbb{D}} \sum_{\zeta_{\mJ(\mathbb{D})}}
    P_{\mR,\rho} \bigl(I \mid \mkD =
    \mathbb{D}, \omega_{\mJ(\mathbb{D})} = \zeta_{\mJ(\mathbb{D})} \bigr) 
    P_{\mR,\rho} \bigl( \mkD = \mathbb{D},
    \omega_{\mJ(\mathbb{D})} = \zeta_{\mJ(\mathbb{D})} \bigr) \\
  &= \sum_{\mathbb{D}} \sum_{\zeta_{\mJ(\mathbb{D})}}
    P_{\mI(\mathbb{D}),f}( A(\zeta_{\mJ(\mathbb{D})}
    \rho^0_{\mathbb{D}}) \circ_r B(\mathbb{D}) )
    P_{\mR,\rho} \bigl( \mkD = \mathbb{D},
    \omega_{\mJ(\mathbb{D})} = \zeta_{\mJ(\mathbb{D})} \bigr).
    \notag
\end{align}
We would like to apply Theorem \ref{T:FKZ2} to the first probability on the right
side of (\ref{E:yesI}), but $\mI(\mathbb{D})$ need not be a circuit-bounded set. 
However, for each connected component $G$ of the interior of
$\mI(\mathbb{D})^{\solid}$,
$\mB(\oG)$ is circuit-bounded, the relevant circuit being the
boundary of $G$.   We let
$\mI_{\main}(\mathbb{D})$ denote the union of all such $\mB(\oG)$. 
Then $\mI_{\main}(\mathbb{D})$ is a finite union of circuit-bounded sets, and
$d_{\mR}(\mI_{\main}(\mathbb{D})^c,\mI_{r/4}(\mathbb{D})) \geq r/4$. Since
$\mI(\mathbb{D})$ has SLC components, for each $\mB(\oG)$, each connected component
of $\mI(\mathbb{D}) \bs \mB(\oG)$ can intersect $\mB(\oG)$ 
in at most a single site. 
This means that under a free boundary condition on
$\mI(\mathbb{D})$, for each $\oG$, regardless of the configuration on
$\mI(\mathbb{D}) \bs
\mB(\oG)$, the effective boundary condition on $\mB(\oG)$ is free.  Therefore
\begin{align} \label{E:mainpart}
  P_{\mI(\mathbb{D}),f}&( A(\zeta_{\mJ(\mathbb{D})}
    \rho^0_{\mathbb{D}}) \circ_r B(\mathbb{D}) ) \\
  &= \sum_{\zeta_{\mI(\mathbb{D}) \bs \mI_{\main}(\mathbb{D})}}
    P_{\mI_{\main}(\mathbb{D}),f} \bigl( A(\zeta_{\mJ(\mathbb{D})}
    \rho^0_{\mathbb{D}}
    \zeta_{\mI(\mathbb{D}) \bs \mI_{\main}(\mathbb{D})}) \circ_r B(\mathbb{D}) 
    \bigr) \notag \\
  &\qquad \qquad \qquad \qquad \cdot P_{\mI(\mathbb{D}),f}(\omega_{\mI(\mathbb{D})
    \bs \mI_{\main}(\mathbb{D})}
    = \zeta_{\mI(\mathbb{D}) \bs \mI_{\main}(\mathbb{D})}) \notag
\end{align} 
By Theorem \ref{T:FKZ2} and Remark \ref{R:SLCC} we have 
\begin{align} \label{E:insidebd}
  P_{\mI_{\main}(\mathbb{D}),f} &\bigl( A(\zeta_{\mJ(\mathbb{D})}
    \rho^0_{\mathbb{D}}
    \zeta_{\mI(\mathbb{D}) \bs \mI_{\main}(\mathbb{D})}) \circ_r B(\mathbb{D}) 
    \bigr)  \\
  &\leq (1 + c_{12}
    e^{-\ep_{6}r}) P_{\mI_{\main}(\mathbb{D}),f}( A(\zeta_{\mJ(\mathbb{D})}
    \rho^0_{\mathbb{D}}
    \zeta_{\mI(\mathbb{D}) \bs \mI_{\main}(\mathbb{D})}))
    P_{\mI_{\main}(\mathbb{D}),f}(B(\mathbb{D})). \notag
\end{align}
Since $P$ has exponential decay of connectivity, it has uniform exponential decay
for the class of all SLC subsets \cite{Al01fin} and therefore by Proposition
\ref{P:splitmix},
\begin{equation} \label{E:Bbound}
  P_{\mI_{\main}(\mathbb{D}),f}(B(\mathbb{D})) \leq (1 + c_{41}e^{-\ep_{26}r})
    P_{\mR,\rho}(B(\mathbb{D})) \leq (1 + c_{37}e^{-\ep_{26}r})
    P_{\mR,\rho}(B).
\end{equation}
Combining (\ref{E:yesI})--(\ref{E:Bbound}) we obtain
\begin{align}
  P_{\mR,\rho}&(I_{AB}) \\
  &\leq (1 + c_{42}e^{-\ep_{27}r})P_{\mR,\rho}(B) \notag \\
  &\qquad \cdot \sum_{\mathbb{D}} \sum_{\zeta_{\mJ(\mathbb{D})}}
    \sum_{\zeta_{\mI(\mathbb{D}) \bs \mI_{\main}(\mathbb{D})}}
    P_{\mI_{\main}(\mathbb{D}),f}( A(\zeta_{\mJ(\mathbb{D})}\rho^0_{\mathbb{D}}
    \zeta_{\mI(\mathbb{D}) \bs \mI_{\main}(\mathbb{D})}) ) \notag \\
  &\qquad \qquad \qquad \qquad \cdot P_{\mR,\rho} \bigl( \mkD = \mathbb{D},,
    \omega_{\mJ(\mathbb{D})} = \zeta_{\mJ(\mathbb{D})}, 
    \omega_{\mI(\mathbb{D}) \bs \mI_{\main}(\mathbb{D})} = 
    \zeta_{\mI(\mathbb{D}) \bs \mI_{\main}(\mathbb{D})} \bigr) \notag \\
  &= (1 + c_{42}e^{-\ep_{27}r})P_{\mR,\rho}(B) P_{\mR,\rho}(A). \notag
\end{align}
This and (\ref{E:noI2}) complete the proof.
\end{proof}

The proof of Theorem \ref{T:Ising} is generally similar to that of Theorem
\ref{T:finitevol}, except that the couplings of the top and bottom layers are
obtained by a different construction.  We will need two lemmas to replace
Proposition \ref{P:strongmix}.  

\begin{lemma} \cite{Al01fin} \label{L:IsingRSMZd}
  Let $\mu^{\beta,h}$ be the Ising model at $(\beta,h)$ on $\ZZd$, with $\beta
< \beta_c(d,h)$.  Let
$\mkL$ be a class of subsets of $\ZZd$ which has the neighborhood
component property.  Suppose that the corresponding FK model has uniform
exponential decay of finite-volume connectivities for the class
$\{ \mB(\Lambda): \Lambda \in \mkL \}$ with wired boundary conditions.  Then
$\mu^{\beta,h}$ has the ratio strong mixing property for the class
$\mkL$ and arbitrary boundary conditions.
\end{lemma}

\begin{lemma} \label{L:IsingRSMZ2}
Let $\mu^{\beta,h}$ be the Ising model at $(\beta,h)$ on $\ZZ$ and let $\mkL$ be
the class of all finite SLC subsets of $\ZZ$ with arbitrary
boundary condition.  Suppose that either (a) $\beta <
\beta_c(2,0)$ and $h = 0$, (b) $\beta > \beta_c(2,0)$ and $h \neq 0$, (c) $\beta <
\beta_c(2,h)$ and the corresponding
FK model has exponential decay of connectivities (in infinite volume), or (d)
$\beta > \beta_c(2,h)$ and the corresponding FK model has exponential decay of dual
connectivities (in infinite volume).  Then $\mu^{\beta,h}$ has the ratio strong
mixing property for the  class $\mkL$.
\end{lemma}
\begin{proof}
Under (c) and (d) this is proved in \cite{Al01fin}.  Suppose (a) holds; then the
Ising model has exponential decay of correlations \cite{ABF}, so by
\ref{E:covarconn} the corresponding FK model has exponential decay of
connectivities, and (a) follows from (c).  Next suppose (b) holds; we
may assume $h>0$.  The Ising model then FKG-dominates the plus phase at
$(\beta,0)$. In the plus phase at $(\beta,0)$, the probability that there is
a path from 0 to $x$ on which all sites $y$ have $\sigma_y = -1$ decays
exponentially in $|x|$ \cite{CCS}, so the same is true at $(\beta,h)$.  It follows
that the corresponding ARC model (an alternate graphical representation of the
Ising model---see
\cite{AlARC}) at $(\beta,h)$ has exponential decay of connectivities (in infinite
volume), which in turn implies that $\mu^{\beta,h}$ has the ratio strong mixing
property for the  class $\mkL$ \cite{Al01fin}.
\end{proof}

It is plausible that the exponential decay assumptions in Lemma
\ref{L:IsingRSMZ2}(c) and (d) are valid for all 
$\beta$ below and above $\beta_c(2,h)$, respectively, but this is not known
rigorously for all $h \neq 0$.

We next construct the coupling that will be used in the $A$ and $B$ pairs for the
Ising model.

\begin{lemma} \label{L:RSMcoupling}
  Suppose $\mu$ is an Ising model on $\ZZd$ having the ratio strong
  mixing property for some class $\mkL$ of finite subsets of $\ZZd$ with arbitrary
  boundary conditions. 
  There exist $c_{i},\ep_i$ as follows.  Suppose $\Delta \subset \Lambda \subset
  \ZZd$ with $\Lambda \in \mkL$, $\alpha \in \{ -1,1 \}^{\partial \Lambda}$,
  $r \geq c_{43} \log |\Delta|$, and
  $\eta_{\Delta} \geq \eta_{\Delta}^{\prime} \in \{-1,1\}^{\Delta}$.
  There exists an FKG coupling $\hat{\mu}_{\Lambda,\alpha}$ of
  $\mu_{\Lambda,\alpha}(\sigma_{\Lambda \bs \Delta} \in \cdot \mid
  \sigma_{\Delta} = \eta_{\Delta})$ and 
  $\mu_{\Lambda,\alpha}(\sigma_{\Lambda \bs \Delta}
  \in \cdot \mid
  \sigma_{\Delta} = \eta_{\Delta}^{\prime})$ with the property that
  \begin{align} \label{E:RSMcouple}
    &\hat{\mu}_{\Lambda,\alpha} \bigl( \bigl\{
      (\sigma_{\Lambda \bs \Delta},\sigma_{\Lambda \bs \Delta}^{\prime}):
      \sigma_x \neq \sigma^{\prime}_x
      \text{ for some } x \in \Lambda \bs \Delta^r(\Lambda) \bigr\} \mid 
      \sigma_{\Lambda \bs \Delta^r(\Lambda)} = \eta_{\Lambda \bs \Delta^r(\Lambda)}
      \bigr) \\
    &\qquad \leq c_{44}e^{-\ep_{28}r} \quad \text{for every } 
      \eta_{\Lambda \bs \Delta^r(\Lambda)} \in 
      \{-1,1\}^{\Lambda \bs \Delta^r(\Lambda)}. \notag
  \end{align}
\end{lemma}

If it were not for the
conditioning on $\sigma_{\Lambda \bs\Delta^r(\Lambda)} = \eta_{\Lambda \bs
\Delta^r(\Lambda)}$, Lemma \ref{L:RSMcoupling} would be a standard result saying
that there is at most an exponentially small probability that the two coupled
configurations, given
$\eta_{\Delta}$ and given $\eta_{\Delta^{\prime}}$,  are unequal far from
$\Delta$.  This standard result leaves open the possibility, though, that there are
a few rare ``hard-to-couple-to'' configurations
$\eta_{\Lambda \bs \Delta}$ which greatly increase the probability of unequal
configurations when they occur, say, in the top layer of the coupled
configuration.  Lemma
\ref{L:RSMcoupling} rules out this possibility. 

Lemma \ref{L:RSMcoupling} extends straightforwardly to any bond or site model
having the ratio strong mixing property.  If the model lacks the FKG property, the
coupling will not be an FKG coupling in general.

We will refer to a coupling of the type guaranteed by Lemma \ref{L:RSMcoupling} as
an \emph{RSM coupling}.

\begin{proof}[Proof of Lemma \ref{L:RSMcoupling}]
Let $\Lambda_{\far} = \Lambda \bs \Delta^r(\Lambda)$ and $\Lambda_{\near} = 
\Delta^r(\Lambda) \bs \Delta$;
we will refer to configurations on $\Lambda_{\far}$ and $\Lambda_{\near}$ as far
and near configurations, respectively.  Define a measure $\nu$ on far
configurations by
\[
  \nu(\zeta_{\Lambda_{\far}}) = \min \bigl(
    \mu_{\Lambda,\alpha}( \sigma_{\Lambda_{\far}}
    = \zeta_{\Lambda_{\far}}
    \mid \sigma_{\Delta} = \eta_{\Delta}), 
    \mu_{\Lambda,\alpha}(\sigma_{\Lambda_{\far}}
    = \zeta_{\Lambda_{\far}}
    \mid \sigma_{\Delta} = \eta_{\Delta}^{\prime}) \bigr), 
\]
and let
\[  \nu_{0} = \nu \left( \{-1,1\}^{\Lambda_{\far}} \right),
\]
\[
  \tau(\cdot) = \frac{\mu_{\Lambda,\alpha}(\sigma_{\Lambda_{\far}} 
  \in \cdot \mid
  \sigma_{\Delta} = \eta_{\Delta}) - \nu(\cdot)}{1 - \nu_{0}}, \qquad 
  \tau^{\prime}(\cdot) = 
  \frac{\mu_{\Lambda,\alpha}(\sigma_{\Lambda_{\far}} \in \cdot \mid
  \sigma_{\Delta} = \eta_{\Delta}^{\prime}) - \nu(\cdot)}{1 - \nu_{0}}.
\]
We may assume  $\nu_{0} < 1$, since otherwise the two measures to be coupled are
identical. Then $\nu/\nu_{0}, \tau$ and $\tau^{\prime}$ are probability measures, 
\begin{equation} \label{E:min0}
  \min \bigl( \tau(\zeta_{\Lambda_{\far}}), 
  \tau^{\prime}(\zeta_{\Lambda_{\far}}) \bigr) 
  = 0 \quad \text{for all } \zeta_{\Lambda_{\far}},
\end{equation}
and by the
ratio strong mixing property, 
\begin{equation} \label{E:maxratio}
  (1 - \nu_{0})\max \left( \frac{\tau(\zeta_{\Lambda_{\far}})}
    {\nu(\zeta_{\Lambda_{\far}})}, 
    \frac{\tau^{\prime}(\zeta_{\Lambda_{\far}})}{\nu(\zeta_{\Lambda_{\far}})} 
    \right) \leq
    c_{45}e^{-\ep_{29}r} \qquad \text{ for all } \zeta_{\Lambda_{\far}}.
\end{equation}
Let $(\chi, \chi^{\prime}), \xi$ and $X$ be independent, with 
$\Prob(X = 1) = \nu_{0}, \Prob(X = 0) = 1 - \nu_{0}$, with
$\xi$ having distribution $\nu/\nu_{0}$, and with
$(\chi, \chi^{\prime})$ having as its distribution an FKG coupling of 
$\tau$ and $\tau^{\prime}$, and with
\[
  \Prob(X = 1) = \nu_{0}, \qquad \Prob(X = 0) = 1 - \nu_{0},
\]
where $\Prob$ denotes the distribution of $(\chi,\chi^{\prime},\xi,X)$.
Note that $\chi \neq \chi^{\prime}$, by (\ref{E:min0}).  Set 
\begin{equation}
  (\sigma_{\Lambda_{\far}},\sigma_{\Lambda_{\far}}^{\prime}) = 
  \begin{cases}
    (\xi,\xi), &\text{ if $X = 1$;} \\
    (\chi,\chi^{\prime}), &\text{ if $X = 0$,}
  \end{cases} \notag
\end{equation}
and let $\hat{\mu}_{\Lambda_{\far}}$ be the distribution of
$(\sigma_{\Lambda_{\far}},\sigma_{\Lambda_{\far}}^{\prime})$.  Then
$\hat{\mu}_{\Lambda_{\far}}$ is an FKG coupling of far configurations, and we
have by (\ref{E:maxratio}), for all
$\zeta_{\Lambda_{\far}}$,
\begin{align} \label{E:farconfigs}
  \hat{\mu}_{\Lambda_{\far}}&(\sigma_{\Lambda_{\far}} \neq
    \sigma_{\Lambda_{\far}}^{\prime} \mid 
    \sigma_{\Lambda_{\far}} = \zeta_{\Lambda_{\far}}) \\
  &\leq \frac{ Pr(\chi = \zeta_{\Lambda_{\far}}, X = 0) }{ Pr(\xi =
    \zeta_{\Lambda_{\far}}, X = 1) }
    \notag \\
  &= (1 - \nu_{0}) \frac{ \tau(\zeta_{\Lambda_{\far}}) }
    { \nu(\zeta_{\Lambda_{\far}}) } \notag \\
  & \leq c_{45}e^{-c_{30}r}. \notag
\end{align}
Now we extend $\hat{\mu}_{\Lambda_{\far}}$ to an
FKG coupling $\hat{\mu}_{\Lambda,\alpha}$ of $\mu_{\Lambda,\alpha}(\sigma_{\Lambda
\bs
\Delta} \in \cdot \mid
\sigma_{\Delta} = \eta_{\Delta})$ and 
$\mu_{\Lambda,\alpha}(\sigma_{\Lambda \bs \Delta}
\in \cdot \mid
\sigma_{\Delta} = \eta_{\Delta}^{\prime})$, by specifying that for each choice of
configurations
$\zeta_{\Lambda_{\far}} \geq \zeta_{\Lambda_{\far}}^{\prime}$, the
distribution of $(\sigma_{\Lambda_{\near}}, 
\sigma_{\Lambda_{\near}}^{\prime})$ given 
$\sigma_{\Lambda_{\far}} = \zeta_{\Lambda_{\far}}, 
\sigma_{\Lambda_{\far}}^{\prime} = \zeta_{\Lambda_{\far}}^{\prime}$ 
is given by an FKG coupling of 
$\mu_{\Lambda,\alpha}(\sigma_{\Lambda_{\near}} \in \cdot \mid
\sigma_{\Delta} = \eta_{\Delta}, \sigma_{\Lambda_{\far}} = 
\zeta_{\Lambda_{\far}})$ and 
$\mu_{\Lambda,\alpha}(\sigma_{\Lambda \bs \Delta}
\in \cdot \mid
\sigma_{\Delta} = \eta_{\Delta}^{\prime}, \sigma_{\Lambda_{\far}} = 
\zeta_{\Lambda_{\far}}^{\prime})$.
It is easily seen that
\begin{align} \label{E:RSMcouple2}
    &\hat{\mu}_{\Lambda,\alpha} \bigl( \bigl\{
      (\sigma_{\Lambda \bs \Delta},\sigma_{\Lambda \bs \Delta}^{\prime}):
      \sigma_x \neq \sigma^{\prime}_x
      \text{ for some } x \in \Lambda_{\far} \bigr\} \mid 
      \sigma_{\Lambda_{\far}} = \eta_{\Lambda_{\far}}
      \bigr) \\
    &\qquad \leq c_{46}e^{-\ep_{30}r} \quad \text{for every } 
      \eta_{\Lambda_{\far}} \in 
      \{-1,1\}^{\Lambda_{\far}}, \notag
\end{align}
since this is only a statement about the coupling of far configurations,
equivalent to (\ref{E:farconfigs}). However, under $\hat{\mu}_{\Lambda,\alpha}$ the
near configuration in
$\sigma_{\Lambda \bs \Delta}$ and the far configuration in $\sigma_{\Lambda \bs
\Delta}^{\prime}$ are conditionally independent given the far configuration in
$\sigma_{\Lambda \bs \Delta}$.  This means that the probability on the left side
of (\ref{E:RSMcouple2}) is unchanged if the conditioning is changed from
$\sigma_{\Lambda_{\far}} = \eta_{\Lambda_{\far}}$ to $\sigma_{\Lambda
\bs \Delta} = \eta_{\Lambda \bs \Delta}$.  Thus (\ref{E:RSMcouple2}) is equivalent
to (\ref{E:RSMcouple}).
\end{proof}

\begin{proof}[Proof of Theorem \ref{T:Ising}.]
We follow the method of (\ref{E:inducstep})--(\ref{E:suffices}), but we use a
different coupling within the $A$ and $B$ pairs.  Recall (see the discussion after
(\ref{E:Markovineq})) that the key property of the coupling for bond models was
that $\omega_{\mS}^1 =  \omega_{\mS}^0$ outside
$C_e = C_e((\zeta_{\mU},1,\omega_{\mS}^1))$, which is the cluster of $e$ for the
top layer.  In the present case, we do not couple with agreement outside a
specified cluster.  Instead, by Lemma \ref{L:IsingRSMZd}, $\mu$ has the ratio
strong mixing property for the class $\mkL$ with arbitrary boundary conditions, so
Lemma \ref{L:RSMcoupling} guarantees the existence of an RSM coupling of the
measures
$\mu_{\Lambda,\eta}( \cdot \mid \sigma_{\Gamma} = \zeta_{\Gamma}, \sigma_x = 1)$
and $\mu_{\Lambda,\eta}( \cdot \mid \sigma_{\Gamma} = \zeta_{\Gamma}, \sigma_x =
0)$, where $\Gamma$ is the set of unsplit sites (the analog of $\mU$) and $x$ is
the site currently being split (the analog of $e$.)  Using an RSM coupling
guarantees that the analog of (\ref{E:suffices}) holds, which, as in the bond
case, leads to (\ref{E:allsplit}), completing the proof of (i).  Then (ii) follows
as in the proof of Theorem \ref{T:FKZd}.
\end{proof}

It should be pointed out that we cannot use an RSM coupling in the proofs of our
theorems on bond models, because typically the ratio strong mixing property will
not apply to the measures $P(\omS \in \cdot \mid \omU = \zeU)$, unless the
configuration $\zeU$ is a special type.  In the FK case, for example, with $P =
P_{\mR,\rho}$ for some
$\mR$ and $\rho$, the ratio strong mixing property cannot be guaranteed unless the
effective boundary condition $(\rho_{\mR^c}\ \zeU)$ on
$(\mS \cup \{e\})^c$ is unique-cluster, which it will not be, for typical $\zeU$.
As mentioned previously, the possible failure of ratio strong mixing is due to the
phenomenon of tunneling, discussed in \cite{Al01fin}.  An RSM coupling can be used
in the Ising case only because the ratio strong mixing property holds for
arbitrary boundary conditions.

\begin{proof}[Proof of Theorem \ref{T:Ising2}.]
By Lemma \ref{L:IsingRSMZ2}, $\mu$ has the ratio strong mixing property for the
class $\mkL$ with arbitrary boundary conditions.  Let $\Theta$ be a
$d_{\mB(\Lambda)}$-ball of radius $\diam_{\mB(\Lambda)}(\Delta)$ centered at some
site in $\Delta$.  Provided $c_{17}$ is large enough, we then have $r \geq c_{15}
\log |\Theta|$, for the $c_{15}$ of Theorem \ref{T:Ising}.
Now the proof can be
completed similarly to that of Theorem \ref{T:Ising}.
\end{proof}


\begin{thebibliography}{99}

\bibitem{ABF}  Aizenman, M., Barsky, D.J., Fernandez, R., \emph{The phase
transition in a general class of Ising-type models is sharp}, J. Stat. Phys.
\textbf{47} (1987), 343--374.

\bibitem{ACCN}  Aizenman, M., Chayes, J.T., Chayes, L., and Newman, 
C.M., \emph{Discontinuity of the magnetization in the} $1/|x - y|^{2}$
\emph{Ising and Potts models}, J. Stat. Phys. \textbf{50} (1988), 1-40.

\bibitem{Al01fin} Alexander, K.S., \emph{Mixing properties and exponential decay
for lattice systems in finite volumes},
www.ma.utexas.edu/mp\_arc-bin/mpa?yn=01-309 (2001).

\bibitem{Al00wu} Alexander, K.S., \emph{Cube-root boundary fluctuations for
droplets in random cluster models}, Commun. Math. Phys. \textbf{224} (2001),
733--781, arXiv:math.PR/0008217 .

\bibitem{AlARC} Alexander, K.S., \emph{The asymmetric random
cluster model and comparison of Ising and Potts models},
Probab. Theory Rel. Fields \textbf{120} (2001), 395--444.

\bibitem{Al00sp} Alexander, K.S., \emph{The spectral gap the the } 2-$D$
\emph{stochastic Ising model with nearly single-spin boundary conditions},
J. Stat. Phys. \textbf{104} (2001), 59--87.

\bibitem{Al98mix} Alexander, K.S., \emph{On weak mixing
in lattice models}, Probab. Theory Rel. Fields \textbf{110}
(1998), 441-471.

\bibitem{AC} Alexander, K.S. and Chayes, L., \emph{Non-perturbative criteria
for Gibbsian uniqueness}, Commun. Math. Phys. \textbf{189} (1997), 447-464.

\bibitem{BBCK} Biskup, M., Borgs, C., Chayes, J. T. and 
Koteck\'y, R., \emph{Gibbs
states of graphical representations of the Potts model with external fields.
Probabilistic techniques in equilibrium and 
nonequilibrium statistical physics.}, 
J. Math. Phys. \textbf{41} (2000), 1170--1210. 

\bibitem{CCS} Chayes, J.T., Chayes, L. and Schonmann, R., \emph{Exponential
decay of connectivities in the two-dimensional Ising model}, J. Stat. Phys.
\textbf{49} (1987), 433-445. 

\bibitem{ES} Edwards, R.G. and Sokal, A.D., \emph{Generalization of the
Fortuin-Kasteleyn-Swendsen-Wang representation and Monte Carlo algorithm},
Phys. Rev. D \textbf{38} (1988), 2009-2012.

\bibitem{Fo1} Fortuin, C.M., \emph{On
the random cluster model. II. The percolation model}, Physica \textbf{58}
(1972), 393-418.

\bibitem{Fo2} Fortuin, C.M., \emph{On
the random cluster model. III. the simple random-cluster process}, Physica
\textbf{59} (1972), 545-570.

\bibitem{FK} Fortuin, C.M. and Kasteleyn, P.W., \emph{On
the random cluster model. I. Introduction and relation to other models},
Physica \textbf{57} (1972), 536-564.

\bibitem{Gr95} Grimmett, G.R., \emph{The stochastic
random-cluster process and uniqueness of random-cluster measures},
Ann. Probab. \textbf{23} (1995), 1461-1510.

\bibitem{Gr95b} Grimmett, G.R., \emph{Comparison and disjoint-occurrence
inequalities for random-cluster models}, J. Stat. Phys. \textbf{78} (1995),
1311-1324.

\bibitem{Gr97}  Grimmett, G.R., \emph{Percolation and disordered
systems}, in \emph{Lectures on Probability Theory and Statistics.
Lectures from the 26th Summer School on Probability Theory held 
in Saint-Flour, August 19--September 4, 1996} (P. Bernard, ed.),
153-300, Lecture Notes in Mathematics \textbf{1665} (1997).

\bibitem{Ha}  Harris, T., \emph{A lower bound for the critical 
probability in a certain percolation process}, Proc. Camb. Phil. Soc.
\textbf{56} (1960), 13-20.

\bibitem{Ho} Holley, R., \emph{Remarks on the FKG inequalities}, Commun. Math.
Phys. \textbf{36} (1974), 227-231.

\bibitem{LMR}  Laanait, L., Messager, A. and Ruiz, J., \emph{Phase
coexistence and surface tensions for the Potts model}, Commun. Math.
Phys. \textbf{105} (1986), 527-545.

\bibitem{Me} Menshikov, M.V., \emph{Coincidence of critical points in percolation
problems}, Soviet Math. Doklady \textbf{33} (1986), 856--859 [Dokl. Acad. Nauk.
SSSR \textbf{288}, 1308--1311, in Russian.] 

\bibitem{MOS} Martinelli, F., Olivieri, E., and Schonmann, R.H., \emph{For 2-D
lattice spin systems weak mixing implies strong mixing}, Commun. Math. Phys.
\textbf{165} (1994), 33-47.

\bibitem{Ne} Newman, C.M., \emph{Disordered Ising systems and random cluster
representations}, in \emph{Probability and Phase Transition (Cambridge 1993)},
247-260, NATO Adv. Sci. Inst. Ser. C Math. Phys. Sci. \textbf{420}, Kluwer
Acad. Publ., Dordrecht (1994). 

\bibitem{Re} Reimer, D., \emph{Proof of the van den Berg--Kesten conjecture},
Combin. Probab. Comput. \textbf{9} (2000), 27-32.

\bibitem{vdBF} van den Berg, J. and Fiebig, U., \emph{On a combinatorial conjecture
concerning disjoint occurrences of events}, Ann. Probab. \textbf{15} (1987),
354--374. 

\bibitem{vdBK} van den Berg, J. and Kesten, H., 
\emph{Inequalities with applications to percolation and
reliability}, J. Appl. Probab. \textbf{22} (1985), 556-569. 

\end{thebibliography}
\end{document}